\documentclass[aos,MSNbibl,nameyear,dvips]{arximspdf}
\usepackage{graphicx}
\doi{10.1214/14-AOS1252} 
\volume{42}
\issue{6}
\pubyear{2014}
\firstpage{2301}
\lastpage{2339}
\docsubty{FLA}

\makeatletter
\newcommand{\rrvert}{\vert}
\newcommand{\rrVert}{\Vert}
\newcommand{\llvert}{\vert}
\newcommand{\llVert}{\Vert}
\newtheorem{theorem}{Theorem}
\newtheorem{lemma}[theorem]{Lemma}
\newproclaim{example}[theorem]{Example}
\newtheorem{corollary}[theorem]{Corollary}
\newproclaim{definition}[theorem]{Definition}
\newproclaim{ass}{Assumption}
\newproclaim{rems}{Remarks}
\newproclaim{rem}{Remark}
\newproclaim{ou}{Outliers}
\newcommand{\R}{\mathbb{R}}
\newcommand{\X}{\mathbb{X}}
\newcommand{\mm}[1]{\ifmmode{#1}\else{\mbox{\(#1\)}}\fi}
\makeatother

\begin{document}
\begin{frontmatter}

\title{Confidence sets for persistence diagrams}
\runtitle{Inference for homology}

\begin{aug}
\author[A]{\fnms{Brittany Terese}~\snm{Fasy}\tsup{*,}\thanksref{T2}\ead[label=e2]{brittany.fasy@alumni.duke.edu}},
\author[B]{\fnms{Fabrizio}~\snm{Lecci}\tsup{$\dagger$,}\thanksref{T3}\ead[label=e3]{lecci@cmu.edu}},
\author[B]{\fnms{Alessandro}~\snm{Rinaldo}\tsup{$\dagger$,}\thanksref{T3}\ead[label=e4]{arinaldo@cmu.edu}},
\author[B]{\fnms{Larry}~\snm{Wasserman}\corref{}\tsup{$\dagger$,}\thanksref{T5}\ead[label=e6]{larry@cmu.edu}},\\
\author[C]{\fnms{Sivaraman} \snm{Balakrishnan}\tsup{$\dagger$}\ead[label=e1]{sbalakri@cs.cmu.edu}}
\and
\author[C]{\fnms{Aarti} \snm{Singh}\tsup{$\dagger$}\ead[label=e5]{aarti@cs.cmu.edu}}
\runauthor{B. T. Fasy et al.}
\affiliation{Tulane University* and Carnegie Mellon University$^{\dagger}$}
\address[A]{B.~T. Fasy\\
Computer Science Department\\
Tulane University\\
New Orleans, Louisiana 70118\\
USA\\
\printead{e2}}
\address[B]{F. Lecci\\
A. Rinaldo\\
L. Wasserman\\
Department of Statistics\\
Carnegie Mellon University\\
Pittsburgh, Pennsylvania 15213\\
USA\\
\printead{e3}\\
\phantom{E-mail: }\printead*{e4}\\
\phantom{E-mail: }\printead*{e6}}
\address[C]{S. Balakrishnan\\
A. Singh\\
Computer Science Department\\
Carnegie Mellon University\\
Pittsburgh, Pennsylvania 15213\\
USA\\
\printead{e1}\\
\phantom{E-mail: }\printead*{e5}}
\end{aug}
\thankstext{T2}{Supported in part by NSF Grant CCF-1065106.}
\thankstext{T3}{Supported in part by NSF CAREER Grant DMS-11-49677.}
\thankstext{T5}{Supported in part by Air Force Grant FA95500910373,
NSF Grant DMS-08-06009.}

\received{\smonth{3} \syear{2013}}
\revised{\smonth{6} \syear{2014}}

%
\begin{abstract}
Persistent homology is a method for
probing topological properties of point clouds
and functions.
The method involves
tracking the birth and death of topological features (2000)
as one varies a tuning parameter.
Features with short lifetimes are informally
considered to be ``topological noise,'' and those with a long lifetime are
considered to be ``topological signal.''
In this paper, we bring some statistical ideas
to persistent homology.
In particular, we derive confidence sets
that allow us to separate topological signal from topological noise.
\end{abstract}

%
\begin{keyword}[class=AMS]
\kwd[Primary ]{62G05}
\kwd{62G20}
\kwd[; secondary ]{62H12}
\end{keyword}
\begin{keyword}
\kwd{Persistent homology}
\kwd{topology}
\kwd{density estimation}
\end{keyword}
\end{frontmatter}

\setcounter{footnote}{3}

\section{Introduction}\label{sectionintro}

Topological data analysis (TDA)
refers to a collection of methods for
finding topological structure in data
[\citet{carlsson2009topology,edelsbrunner2010computational}].
TDA has been used in protein analysis,
image processing, text analysis, astronomy, chemistry and
computer vision, as well as in other fields.

%
\begin{figure}[t]

\includegraphics{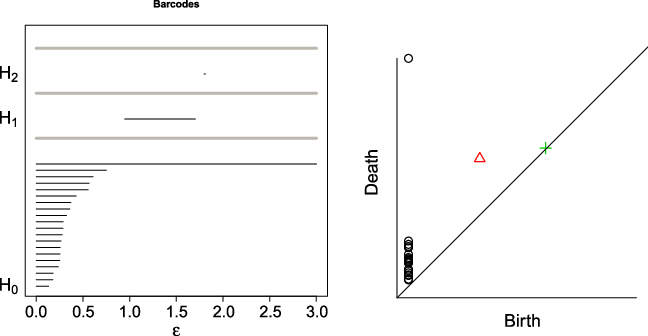}

\caption{(Left) barcode plot. Each horizontal bar shows the lifetime of a
topological feature as a function of $\varepsilon$
for the homology groups $H_0$, $H_1$, and $H_2$. Significant features
have long horizontal bars.
(Right) persistence diagram. The points in the persistence diagram are in
one-to-one correspondence with the bars in the barcode plot. The birth and death
times of
the
barcode become the $x$- and $y$-coordinates of the persistence diagram.
Significant features are far from the diagonal.}\label{figBarcode}
\end{figure}

%
\begin{figure}[b]

\includegraphics{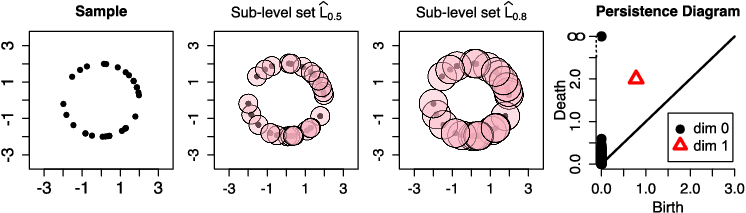}

\caption{(Left) 30 data points $\mathcal{S}_{30}$ sampled from the
circle of radius 2. (Middle left) sub-levels
set \mbox{$\widehat L_{0.5}=\{x\dvtx d_{\mathcal{S}_{30}}\leq0.5\}$}; at
$\varepsilon=0.5$, the sub-level set consists of two connected
components and zero loops. (Middle right) sub-levels set $\widehat
L_{0.8}=\{x\dvtx
d_{\mathcal{S}_{30}} \leq0.8\}$; as we keep increasing $\varepsilon
$, we assist
at the birth and death of topological features; at $\varepsilon=0.8$,
one of the
connected components dies (is merged with the other one), and a one-dimensional
hole appears; this loop will die at $\varepsilon=2$, when the union of
the pink
balls
representing the distance function becomes simply connected. Right: the
persistence diagram summarizes the topological features of
the sampled points. The black dots represent the connected components: 30
connected components are present at $\varepsilon=0$, and they
progressively die
as $\varepsilon$ increases, leaving only one connected component. The red
triangle represents
the unique one-dimensional hole that appears at $\varepsilon=0.8$ and
dies at
$\varepsilon=2$.}\label{figsublevel}
\end{figure}

One approach to TDA is \emph{persistent homology}, a branch of
computational topology that leads to a plot called a \emph{persistence diagram}.
This diagram can be thought of as a summary statistic, capturing
multi-scale topological
features. This paper studies the statistical properties of persistent homology.
Homology detects the connected components,
tunnels, voids, etc., of a topological space ${\mathbb M}$. Persistent
homology measures these
features by assigning a birth and a death value to each feature.
For example, suppose we sample $\mathcal{S}_n=\{X_1,\ldots,X_n\}$ from
an unknown distribution $P$. We are
interested in estimating the homology of the support of $P$. One method for
doing so would be to take the union of the set of balls centered at the points
in $\mathcal{S}_n$ as we do in the supplementary material [\citet{supplement}]; however, to do\vadjust{\goodbreak}
so, we must choose a radius for these balls. Rather than
choosing a single radius, we use persistent homology to summarize the
homology for all radii, as shown in Figures~\ref{figBarcode}~and~\ref{figsublevel}.
Persistence
points far from the diagonal $x=y$ represent topological features
common to a
long interval of radii.
We define persistent homology precisely in Section~\ref{sectionreview}.


One of the key challenges in persistent homology
is to find a way to separate the noise from the signal in the
persistence diagram, and this paper suggests several statistical
methods for
doing so. In particular, we
provide a confidence set for the
persistence diagram ${\mathcal P}$
corresponding to the distance function to a topological space ${\mathbb M}$
using a sample from a distribution $P$ supported on
${\mathbb M}$.
The confidence set has the form
${\mathcal C}=\{ {\mathcal Q}\dvtx  W_\infty(\widehat{\mathcal
P},{\mathcal Q})\leq c_n\}$,
where\vspace*{1pt}
${\mathcal Q}$ varies over the set of all persistence diagrams,
$\widehat{\mathcal P}$ is an estimate of the persistence
diagram constructed from the sample,
$W_\infty$ is a metric on the space of
persistence diagrams, called the bottleneck distance,
and $c_n$ is an appropriate, data-dependent quantity.
In addition, we study the upper
level sets of density functions by computing confidence sets for their
persistence~diagrams.

\subsection*{Goals}
There are two main goals in this paper.
The first is to introduce persistent homology to statisticians.
The second is to derive confidence sets for certain key quantities in
persistent homology.
In particular, we derive a simple method
for separating topological noise from topological signal.
The method has a simple visualization: we only need to add a band
around the
diagonal of the persistence diagram. Points in the band are consistent with
being noise.
We focus on simple, synthetic examples in this
paper as proof of concept.

\subsection*{Related work}
Some key references
in computational topology
are
[\citet{bubenik2007statistical,carlsson2009topology,ghrist2008barcodes,carlsson2009theory,edelsbrunner2008persistent,chazal2008towards,chazal2011persistence}].
An \hyperref[sectionintro]{Introduction} to homology can be found in
\citet{hatcher2002algebraic}.
The probabilistic basis for
random persistence diagrams is studied in
\citet{mileyko2011probability} and
\citet{turner2012fr}.
Other relevant probabilistic results can be found in
\citeauthor{kahle2009topology} (\citeyear{kahle2009topology,kahle2011random}), \citet{penrose2003random,kahle2010limit}.
Some statistical results for homology and persistent homology include
\citet{bubenik2010statistical,blumberg2012persistent,balakrishnan2011minimax,joshi2010comparing,bendich2010towards,chazal2010geometric,bendich2011improving},
\citeauthor{niyogi2008finding} (\citeyear{niyogi2008finding,niyogi2011topological}).
The latter paper considers a challenging example
which involves a set of the form in Figure~\ref{figGlassesUniform}
(top left), which we consider later in the paper.
\citet{heo2012topological}
contains a detailed analysis
of data on a medical procedure for upper jaw expansion
that uses persistent homology as part of a
nonlinear dimension reduction of the data.
\citet{chazal2013optimal} is closely related to this paper;
they find convergence rates for persistence diagrams computed from data sampled
from a
distribution on a metric space. As the authors point out, finding
confidence intervals that follow from their methods is a challenging problem
because their probability inequalities involve unknown constants.
Restricting our attention to manifolds embedded in $\R^D$, we are able to
provide several methods to compute confidence sets for persistence diagrams.

\subsection*{Outline}
We define persistent homology formally in Section~\ref{sectionreview} and
provide additional details in the supplementary material [\citet{supplement}].
The statistical model is defined in Section~\ref{sectionmodel}.
Several methods for constructing confidence intervals
are presented in Section~\ref{sectionconfidence}.
Section~\ref{sectionexperiments} illustrates the ideas
with a few numerical experiments.
Proofs are contained in Section~\ref{sectionproofs}.
Finally,
Section~\ref{sectionconclusion} contains concluding remarks.

\subsection*{Notation}
We write $a_n \preceq b_n$ if there exists $c>0$ such that
$a_n \leq c b_n$ for all large $n$. We write $a_n \asymp b_n$ if $a_n
\preceq b_n$ and $b_n \preceq a_n$.
For any $x\in\mathbb{R}^D$ and any $r\geq0$,
$B(x,r)$ denotes the $D$-dimensional ball of radius $r>0$ centered at $x$.
For any closed set $A \subset\R^D$, we define
\emph{the} (\emph{Euclidean}) \emph{distance~function}
%
\begin{equation}
\label{eqdist} d_A(x) = \inf_{y\in A} \llVert y-x
\rrVert _2.
\end{equation}
In addition, for any $\varepsilon\geq0$, the Minkowski sum is defined as
%
\begin{equation}
A\oplus\varepsilon= \bigcup_{x\in A}B(x,\varepsilon) =
\bigl\{ x\dvtx  d_A(x) \leq \varepsilon \bigr\}.
\end{equation}
The \emph{reach}\label{labelreach}
of $A$---denoted by $\operatorname{reach}(A)$---is the largest $\varepsilon\geq0$
such that each point in $A\oplus\varepsilon$ has a unique
projection onto $A$ [\citet{federer}].
If $f$ is a real-valued function, we define the
\emph{upper level set} $\{x\dvtx  f(x) \geq t\}$,
\emph{the lower level set} $\{x\dvtx  f(x) \leq t\}$ and
\emph{the level set} $\{x\dvtx  f(x) =t \}$.
If
$A$ is measurable, we write
$P(A)$ for the probability of $A$.
For more general events $A$,
$\mathbb{P}(A)$ denotes
the probability of $A$ on an appropriate probability space.
In particular, if $A$ is an event in the $n$-fold probability space under
random sampling, then
$\mathbb{P}(A)$ means probability under
the product measure
$P\times\cdots\times P$. In some places, we use symbols like
$c,c_1,C,\ldots,$ as generic positive constants.
Finally, if two sets $A$ and $B$ are homotopic, we write
$A\cong B$.

\section{Brief introduction to persistent homology}\label{sectionreview}

In this section, we provide a brief
overview of persistent homology.
In the supplementary material [\citet{supplement}], we provide a more details on relevant concepts from
computational topology; however,
for a more complete coverage of persistent homology, we refer the reader
to~\citet{edelsbrunner2010computational}.

Given a real-valued function $f$, persistent homology describes
how the topology of the lower level sets $f^{-1}(- \infty, t]$
change as $t$ increases from $-\infty$ to $\infty$. In particular, persistent
homology describes $f$ with a multiset
of points in the plane, each corresponding to
the birth and death of a homological
feature that existed for some interval of $t$.

First, we consider the case where $f$
is a distance function.
Let $K$ be a compact subset of $\mathbb{R}^D$,
and let $d_K \dvtx \R^D \to\R$ be the distance function to $K$.
Consider the sub-level set
$L_t = \{x\dvtx  d_K(x) \leq t\}$; note that $K=L_0$.
As $t$ varies from 0 to~$\infty$,
the set $L_t$ changes.
Persistent homology summarizes how the topological
features of~$L_t$ change as a function of $t$.
Key topological features of a set include
the connected components (the zeroth order homology),
the tunnels (the first order homology),
voids (second order homology), etc.
These features can appear (be born) and disappear (die) as $t$ increases.
For example, connected components of $L_t$
die when they merge with other connected~components.

Each topological feature
has a birth time $b$ and a death time $d$.
In general, there will be
a set of features
with birth and death times
$(b_1,d_1),\ldots, (b_m,d_m)$.
These points can be plotted on the plane,
resulting in a persistence diagram ${\mathcal P}$; see Figures~\ref{figBarcode} and~\ref{figsublevel}.
Alternatively, we can represent the
pair $(b_i,d_i)$ as the interval $[b_i,d_i]$. The set of intervals is referred
to as a \emph{barcode plot}; see Figure~\ref{figBarcode}. We view the
persistence diagram and the barcode plot as topological summaries of
the input
function or data.
Points near the diagonal in the persistence diagram (i.e., the short
intervals in the barcode plot) have short lifetimes
and are considered ``topological noise.'' Most applications are
interested in
features that we can distinguish from noise; that is, those features that
persist for a large range of values $t$.

\subsection{Persistent homology}
We present persistent homology as a summary of the input function
or data, as the goal of this paper is to define methods for computing
confidence sets for that summary.

Given data points ${\mathcal S}_n = \{ X_1, \ldots, X_n \}$, we are interested
in understanding the homology of the $d$-dimensional compact
topological space
${\mathbb M}\subset\R^D$ from which the data were sampled; see
the supplementary material [\citet{supplement}] for the
definition of homology. If our sample is dense
enough, and the topological space has a nice embedding in~$\R^D$,
then $H_p({\mathbb M})$ is a\vspace*{2pt} subgroup of the $p$th homology group of
the sublevel set $\widehat{L}_{\varepsilon} = \{x\dvtx  d_{{\mathcal
S}_n}(x) \leq\varepsilon\} $
for an interval of values of $\varepsilon$.
Choosing the right
$\varepsilon$ is a difficult task: small $\varepsilon$
will have the homology of $n$ points, and large $\varepsilon$
will have the homology of a single point. Using persistent homology,
we avoid choosing a single $\varepsilon$ by assigning a persistence
value to
each nontrivial
homology generator that is realized as $\widehat{L}_{\varepsilon}$
for some
nonnegative $\varepsilon$.
This persistence value is defined to be the length of the interval of
$\varepsilon$ for which that feature occurs. See Figure~\ref{figsublevel}.

To consider $\widehat{L}_{\varepsilon}$
for every $\varepsilon$ in $(0, \infty)$ would be infeasible. Hence, we
restrict our attention to equivalence classes of homology groups.
Since $H(\widehat{L}_{r}) = H({\rm\check{C}ech}({\mathcal S}_n,
r))$, we use the
\v{C}ech complex to compute the homology of the lower level sets; see
the supplementary material [\citet{supplement}] for the definition of a \v{C}ech complex.
Let $r_1, \ldots, r_k$ be the set of radii such that the complexes
${\rm\check{C}ech}({\mathcal S}_n, r_i)$ and ${\rm\check{C}ech}({\mathcal S}_n,
r_i-\varepsilon)$ are not identical
for sufficiently small~$\varepsilon$. Letting\vspace*{1pt} $K_0 = \varnothing$,
$K_{k+1}$ be the maximal simplicial complex defined on ${\mathcal
S}_n$, and
$K_i = {\rm\check{C}ech}({\mathcal S}_n, (r_i+r_{i-1})/2)$, the
sequence of
complexes is the \emph{\v{C}ech filtration} of~$d_{\mathcal{S}_n}$.
For all
$s < t$, there is a natural inclusion \mbox{$i_{s,t} \dvtx  K_s
\hookrightarrow K_t$}
that induces a group homomorphism $i_{s,t}^* \dvtx  H_p(K_s) \to
H_p(K_t)$. Thus we have the following sequence of homology groups:
%
\begin{equation}
\label{eqinducedhomology} H_p\bigl(\llvert K_0\rrvert \bigr) \to
H_p\bigl(\llvert K_1\rrvert \bigr) \to\cdots\to
H_p\bigl(\llvert K_n\rrvert \bigr).
\end{equation}
We say that a homology class $[\alpha]$ represented by a $p$-cycle
$\alpha$ is
\emph{born} at $K_s$ if $[\alpha]$ is not supported in~$K_r$
for any $r < s$, but is nontrivial in $H_p(\llvert K_s\rrvert )$.
The class $[\alpha]$ born at~$K_s$
\emph{dies} going into $K_t$ if $t$ is the smallest index such that
the class $[\alpha]$ is supported in the
image of~$i_{s-1,t}^*$.
The birth at $s$ and death at $t$ of $[\alpha]$ is
recorded as the point~$(s,t)$ in the $p$th persistence diagram
$\mm{\mathcal{P}_{{p}}({d_{\mathcal{S}_n}})}$,
which we now formally define.

\begin{figure}

\includegraphics{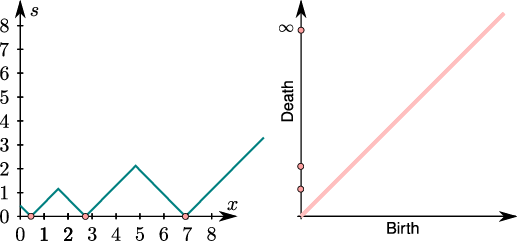}

\caption{(Left) distance function $d_{\X}$ for a set of three points.
(Right) corresponding persistence diagram.
If $\mathbb{X} \subset\R^1$, then the persistence
diagram records the death times of the initial components. Since
all components are born at $s=0$, all points in the persistence
diagram appear on the $y$-axis. The point labeled with $\infty$
represents an essential class, that is, one that is inherited
from the domain (in this case, $\R^1$).}\label{figpersR1}
\end{figure}

%
\begin{definition}[(Persistence diagram)]
Given a function $f \dvtx \X\to R$, defined for a triangulable
subspace of
$\R^D$,
the $p$th \emph{persistence diagram} $\mm{\mathcal{P}_{{p}}({f})}$
is the multiset of points in the extended plane $\overline{\R}^2$, where
$\overline{\R}=\R\cup\{-\infty, +\infty\}$, such that the each point
$(s,t)$ in the diagram represents a distinct
topological feature that existed in $H_p({f}^{-1}((-\infty,r]))$ for
$r \in[s,t)$. The
\emph{persistence barcode} is a multiset of intervals that
encodes the same information as the persistence
diagram by representing the point $(s,t)$ as the interval~$[s,t]$.
\end{definition}
In sum, the zero-dimensional diagram
$\mm{\mathcal{P}_{{0}}({f})}$ records the birth and death of
components of the lower
level
sets; more generally, the $p$-dimensional diagram
$\mm{\mathcal{P}_{{p}}({f})}$ records the $p$-dimensional holes of
the lower level
sets. We let $\mm{\mathcal{P}({f})}$ be the overlay of all
persistence diagrams for $f$;
see Figures~\ref{figsublevel}~and~\ref{figpersR1},
for examples, of persistent homology of one-dimensional
and two-dimensional distance~functions.

%

In the supplementary material [\citet{supplement}], we see that the homology of a lower
level set is equivalent to a ${\rm\check{C}ech}$ complex and can be estimated
by a Vietoris--Rips complex. Therefore, in Section~\ref{sectionexperiments},
we use the Vietoris--Rips filtration to compute the confidence sets for
persistence diagrams of distance functions.

\subsection{Stability}

We say that the persistence diagram is stable if
a small change in the input function produces a
small change in the persistence diagram. There are
many variants of the stability result for persistence diagrams,
as we may define different ways of measuring distance
between functions or distance between persistence diagrams.
We are interested in using the $L_{\infty}$-distance between
functions and the bottleneck distance between persistence diagrams;
see the supplementary material [\citet{supplement}] for the definition of bottleneck distance.

The $L_{\infty}$-distance $\llVert  f-g\rrVert  _{\infty}$ is the maximum
difference between
the
function values
\[
\llVert f-g\rrVert _\infty= \sup_{x \in\X} \bigl\llvert
f(x)-g(x) \bigr\rrvert. %
\]
We can upper bound the
bottleneck distance
between two persistence diagrams by the $L_{\infty}$-distance between the
corresponding functions:
%
\begin{theorem}[(Bottleneck stability)]\label{thmbotStab}
Let $\X$ be finitely triangulable, and let $f,g \dvtx \X\to\R$ be
continuous. Then the bottleneck
distance between the corresponding persistence diagrams is
bounded from above by the $L_{\infty}$-distance between~them,
%
\begin{equation}
\label{eqbotStab} W_\infty\bigl(\mm{\mathcal{P}_({f})}, \mm{
\mathcal{P}_({g})}\bigr) \leq\llVert f-g\rrVert _\infty.
\end{equation}
\end{theorem}

The bottleneck stability theorem is one of the main requirements
for our methods to work.
We refer the reader to \citet{cohenSteiner2007stability} and
to \citet{chazal2012structure} for proofs of this
theorem. Alternatively, one can consider different distance functions.
The Wasserstein distance is defined by
finding the perfect pairing that minimizes the sum (rather than the supremum)
of the pairwise distances. For a restricted class of functions, there
exists a
stability theorem for the Wasserstein distance; see~\citet
{cs2010lipschitz}.
We note that the techniques to compute confidence sets presented
in this paper can be extended to define a confidence set for the persistence
diagram under the Wasserstein distance using a stronger set of assumptions
on~${\mathbb M}$.

\subsection{Hausdorff distance}

Let $A,B$ be compact subsets of $\R^D$. One way to measure
the distance between these sets is to take the \emph{Hausdorff distance},
denoted by $\mm{\mathbf{H}}({A},{B})$, which is the
maximum Euclidean distance from
a point in one set to the closest point in the other set,
\begin{eqnarray*}
\mm{\mathbf{H}}({A},{B}) &=& \max \Bigl\{ \max_{x\in A} \min
_{y\in B} \llVert x-y\rrVert, \max_{x\in B} \min
_{y\in A} \llVert x-y\rrVert \Bigr\}
\\
&=& \inf \{ \epsilon\dvtx  A \subset B\oplus\epsilon\mbox{ and } B \subset A\oplus
\epsilon \}.
\end{eqnarray*}

The stability result
Theorem~\ref{thmbotStab}
is key.
Let ${\mathbb M}$ be a $d$-manifold embedded in a compact subset $\X$
of $\R^D$.
If $S$ is any subset of ${\mathbb M}$,
${\mathcal P}_S$ is the persistence diagram based on the
lower level sets $\{x\dvtx  d_S(x) \leq t\}$, and
${\mathcal P}$ is the persistence diagram based on
the lower level sets $\{x\dvtx  d_{{\mathbb M}}(x) \leq t\}$, then, by
Theorem~\ref{thmbotStab},
%
\begin{equation}
W_\infty({\mathcal P}_S,{\mathcal P}) \leq\llVert
d_S - d_{{\mathbb
M}}\rrVert _\infty= \mm{
\mathbf{H}}({S},{{\mathbb M}}).
\end{equation}
We bound $\mm{\mathbf{H}}({S},{{\mathbb M}})$ to obtain a bound on
$W_\infty
({\mathcal P}_S,{\mathcal P})$.
In particular, we obtain a confidence set for
$W_\infty({\mathcal P}_S,{\mathcal P})$ by
deriving a confidence set for $\mm{\mathbf{H}}({S},{{\mathbb M}})$.
Thus the stability theorem reduces the problem of inferring
persistent homology to the problem of inferring Hausdorff distance.
Indeed, much of this paper is devoted
to the latter problem.
We would like to point out that
the Hausdorff distance plays
an important role
in many statistical problems.
Examples include
\citet{cuevas2009set,cuevas2001cluster},
\citeauthor{cuevas1997plug} (\citeyear{cuevas1997plug,cuevas1998visual}),
\citet{cuevas2012statistical,cuevas2004boundary,mammen1995asymptotical}.
Our methods could potentially be useful
for these problems as well.

\section{Statistical model}\label{sectionmodel}

As mentioned above, we want to estimate the homology of a set ${\mathbb M}$.
We do not observe ${\mathbb M}$ directly;
rather, we observe a sample
${\mathcal S}_n = \{ X_1,\ldots, X_n \}$ from a distribution $P$
that is concentrated on or near ${\mathbb M}\subset\R^D$.
For example, suppose ${\mathbb M}$ is a circle.
Then the homology of the data set
${\mathcal S}_n$ is not equal to the homology
of ${\mathbb M}$; however, the set
$\widehat L_\varepsilon= \{ x\dvtx  d_{{\mathcal S}_n}(x) \leq\varepsilon
\}=\bigcup_{i=1}^n B(X_i,\varepsilon)$,
where $B(x,\varepsilon)$ denotes the Euclidean ball of radius
$\varepsilon$
centered at $x$, captures the homology of ${\mathbb M}$ for an interval of
values $\varepsilon$. Figure~\ref{figsublevel} shows
$\widehat L_\varepsilon$ for increasing values of $\varepsilon$.

Let $X_1,\ldots, X_n \stackrel{\mathrm{i.i.d.}}{\sim} P$
where $X_i \in\mathbb{R}^D$.
Let ${\mathbb M}$ denote the $d$-dimensional support of $P$.
Let us define the following quantities:
%
\begin{equation}
\label{eqrho-t} \rho(x,t) = \frac{P(B(x,t/2))}{t^d}, \qquad \rho(t) = \inf
_{x
\in{\mathbb M}} \rho(x,t).
\end{equation}
We assume that $\rho(x,t)$ is a continuous function of $t$,
and we define
%
\begin{equation}
\label{eqrho} \rho(x, \downarrow0) = \lim_{t\to0} \rho(x,t),
\qquad\rho= \lim_{t\to0} \rho(t).
\end{equation}

Note that if $P$ has continuous density $p$
with respect to the uniform measure on~${\mathbb M}$,
then $\rho\propto\inf_{x\in{\mathbb M}}p(x)$.
Until Section~\ref{sectiondensity},
we make the following assumptions:
\renewcommand{\theass}{A\arabic{ass}}
\begin{ass}\label{assA1}
${\mathbb M}$ is $d$-dimensional compact manifold
(with no boundary), embedded in $\mathbb{R}^D$
and~\mbox{$\operatorname{reach}({\mathbb M})>0$}. (The definition of reach was given
in Section~\ref{sectionintro}.)
\end{ass}

\begin{ass}\label{assA2}
For each $x\in{\mathbb M}$,
$\rho(x,t)$ is a bounded continuous function of $t$,
differentiable for $t\in(0,t_0)$
and right differentiable at zero.
Moreover,
$\partial\rho(x,t)/\partial t$ exists and is bounded away from zero and
infinity
for $t$ in an open neighborhood of zero.
Also, for some $t_0>0$ and some $C_1$ and $C_2$, we have
%
\begin{equation}
\sup_x \sup_{0 \leq t \leq t_0} \biggl\llvert
\frac{\partial\rho
(x,t)}{\partial t} \biggr\rrvert \leq C_1 < \infty\quad\mbox{and}\quad \sup
_{0 \leq t \leq t_0} \bigl\llvert \rho'(t) \bigr\rrvert \leq
C_2 < \infty.
\end{equation}
\end{ass}

\begin{rems*}
The reach of ${\mathbb M}$ does not appear explicitly in the results
as the dependence is implicit and does not affect the rates
in the asymptotics.
Note that
if $P$ has a density $p$ with respect to the
Hausdorff measure on ${\mathbb M}$,
then Assumption~\ref{assA2} is satisfied as long as $p$ is smooth and bounded away from zero.
Assumption~\ref{assA1} guarantees that
as $\varepsilon\to0$, the covering number $N(\varepsilon)$
satisfies
$N(\varepsilon)\asymp(1/\varepsilon)^d$.
However, the conditions are likely stronger than needed.
For example, it suffices that
${\mathbb M}$ be compact and $d$-rectifiable.
See, for example,
\citet{mattila1999geometry} and
\citet{ambrosio2000functions}.

Recall that the distance function is
$d_{{\mathbb M}}(x) = \inf_{y\in{\mathbb M}}\llVert  x-y\rrVert  $, and
let ${\mathcal P}$ be the persistence diagram
defined by the lower level sets
$\{x\dvtx  d_{{\mathbb M}}(x) \leq\varepsilon\}$.
Our target of inference is ${\mathcal P}$.
Let $\widehat{\mathcal P}$ denote the
persistence diagram of the
$\{x\dvtx  d_{{\mathcal S}_n}(x) \leq\varepsilon\}$
where ${\mathcal S}_n= \{X_1,\ldots, X_n\}$.
We regard $\widehat{\mathcal P}$
as an estimate of ${\mathcal P}$.
Our main goal is to find
a confidence interval for
$W_\infty({\mathcal P},\widehat{\mathcal P})$
as this implies a confidence set
for the persistence diagram.

Until Section~\ref{sectiondensity},
we assume that the dimension of ${\mathbb M}$ is known
and that the support of the distribution is ${\mathbb M}$
which is sometimes referred to as the noiseless case.
These assumptions may seem unrealistic to statisticians
but are, in fact, common in computational geometry.
In Section~\ref{sectiondensity},
we weaken the assumptions.
Specifically, we allow outliers, which means there may be points not
on~${\mathbb M}$.
\citet{bendich2011improving} show that methods based on the
\v{C}ech complex perform poorly when there are outliers.
Instead, we estimate the persistent homology of the upper level sets
of the density function.
We shall see that the methods in
Section~\ref{sectiondensity}
are quite robust.
\end{rems*}

\section{Confidence sets}\label{sectionconfidence}

Given $\alpha\in(0,1)$, we will find $c_n \equiv c_n(X_1,\ldots,X_n)$ such
that
%
\begin{equation}
\limsup_{n\to\infty}\mathbb{P}\bigl( W_\infty(\widehat{
\mathcal P},{\mathcal P}) > c_n\bigr) \leq \alpha.
\end{equation}
It then follows that
$C_n = [0,c_n]$ is an asymptotic $1-\alpha$ confidence set for
the bottleneck distance
$W_\infty(\widehat{\mathcal P},{\mathcal P})$, that is,
%
\begin{equation}
\liminf_{n\to\infty}\mathbb{P} \bigl( W_\infty(\widehat{
\mathcal P},{\mathcal P}) \in [0,c_n] \bigr) \geq1-\alpha.
\end{equation}

Recall that,
from Theorem~\ref{thmbotStab} and the fact
that $\llVert  d_{{\mathbb M}}-d_{{\mathcal S}_n}\rrVert  _\infty=\mm
{\mathbf
{H}}({{\mathcal S}_n},{{\mathbb M}})$, we have
%
\begin{equation}
W_\infty(\widehat{\mathcal P},{\mathcal P}) \leq\mm{\mathbf {H}}({{
\mathcal S}_n},{{\mathbb M}}),
\end{equation}
where
${\mathcal S}_n = \{X_1,\ldots, X_n\}$ is the sample
and $\mathbf{H}$ is the Hausdorff distance.
Hence it suffices to find $c_n$ such that
%
\begin{equation}
\limsup_{n\to\infty}\mathbb{P} \bigl(\mm{\mathbf{H}}({{\mathcal
S}_n},{{\mathbb M}})> c_n \bigr) \leq \alpha.
\end{equation}

\begin{figure}

\includegraphics{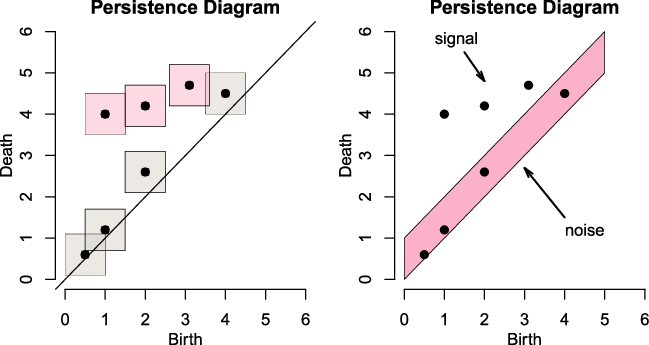}

\caption{First, we obtain the confidence interval $[0,c_n]$ for
$W_\infty(\widehat{\mathcal P},{\mathcal P})$.
If a box of side length $2c_n$ around a point in the diagram hits the diagonal,
we consider that point to be noise.
By putting a band of width $\sqrt{2}c_n$ around the diagonal,
we need only check which points fall inside the band and
outside the band. The plots show the two different ways
to represent the confidence interval $[0, c_n]$. For this particular example
$c_n=0.5$.}\label{figstylized}
\end{figure}

The confidence set ${\mathcal C}_n$ is a subset of all persistence
diagrams whose distance to~$\widehat{\mathcal P}$ is at most $c_n$,
\[
\mathcal{C}_n = \bigl\{ {\widetilde{\mathcal P}}\dvtx W_\infty(\widehat{\mathcal P},{\widetilde{\mathcal P}}) \leq
c_n\bigr\}.
\]
We can visualize ${\mathcal C}_n$ by centering a box of side length $2
c_n$ at each point $p$ on the persistence diagram. The point $p$ is
considered indistinguishable from noise if the corresponding box,
formally defined as $\{q\in\mathbb{R}^2\dvtx  d_\infty(p,q) \leq c_n\}$,
intersects the~diagonal. Alternatively, we can visualize the
confidence set by adding a band of width $\sqrt{2}c_n$ around the
diagonal of the persistence diagram $\widehat{\mathcal P}$.
The interpretation is this: points in the band are not significantly
different from noise. Points above the band can be interpreted as
representing a significant topological feature.
That is, if the confidence set for a point on the diagram hits the diagonal,
then we cannot rule out that the lifetime of that feature is 0, and we
consider it to be noise.
(This is like saying that if a confidence interval for a treatment
effect includes 0,
then the effect is not distinguishable form ``no effect.'')
This leads to the diagrams shown in Figure~\ref{figstylized}.

%

\begin{rem*}
This simple dichotomy of ``signal'' and
``noise'' is not the only way to
quantify the uncertainty in the persistence diagram.
Indeed, some points near the diagonal may represent
interesting structure.
One can imagine endowing each point in the diagram
with a confidence set, possibly of different sizes and shapes.
But for the purposes of this paper,
we focus on the simple method described above.

The first three methods that we present
are based on the persistence diagram constructed from
the \v{C}ech complex.
The fourth method takes a different approach completely
and is based on density estimation.
We define the methods in this section;
we illustrate them in Section~\ref{sectionexperiments}.
\end{rem*}

\subsection{Method \textup{I}: Subsampling}

The first method uses subsampling.
The usual approach to subsampling
[see, e.g., \citet{politis1999subsampling,RomanoShaikh12}]
is based on the assumption that
we have an estimator $\hat\theta$ of a parameter $\theta$ such that
$n^\xi(\hat\theta- \theta)$
converges in distribution to some fixed distribution~$J$
for some $\xi>0$.
Unfortunately,
our problem is not of this form.
Nonetheless, we can still use subsampling
as long as we are willing to have conservative confidence intervals.

Let $b = b_n$ be such that
$b = o ( \frac{n}{\log n}  )$ and $b_n \to\infty$.
We draw all $N$ subsamples
${\mathcal S}_{b,n}^1,\ldots, {\mathcal S}_{b,n}^N$, each of size~$b$,
from the data
where $N={n\choose b}$.
(In practice, as is always the case with subsampling,
it suffices to draw a large number of subsamples randomly
rather than use all $N$ subsamples. But, for the theory, we
assume all $N$ subsamples are used.)
Let
$T_j = \mm{\mathbf{H}}({{\mathcal S}_{b,n}^j},{{\mathcal S}_n})$,
$j=1,\ldots, N$.
Define
%
\begin{equation}
\label{eqLb1} L_b(t) = \frac{1}{N}\sum
_{j=1}^N I(T_j > t),
\end{equation}
and let
$c_b=2 L_b^{-1}(\alpha)$.
Recalling the definition of $\rho$ from~(\ref{eqrho}), we can prove the
following theorem:

%
\begin{theorem}
\label{theoremsubsample}
Assume that $\rho>0$.
Then, for all large $n$,
%
\begin{equation}
\mathbb{P}\bigl( W_\infty(\widehat{\mathcal P},{\mathcal P}) >
c_b\bigr) \leq \mathbb{P} \bigl( \mm{\mathbf{H}}({{\mathcal
S}_n},{{\mathbb M}}) > c_b \bigr) \leq \alpha+ O \biggl(
\frac{b }{n} \biggr)^{1/4}. 
\end{equation}
\end{theorem}

\subsection{Method \textup{II}: Concentration of measure}

The following lemma is similar to theorems
in \citet{devroye1980detection},
\citet{cuevas2001cluster} and
\citet{niyogi2008finding}.

%
\begin{lemma}
\label{lemmaconcentration}
For all $ t>0$,
%
\begin{equation}
\label{eqlam1} \mathbb{P}\bigl( W_\infty(\widehat{\mathcal P},{\mathcal
P}) > t\bigr) \leq \mathbb{P} \bigl(\mm{\mathbf{H}}({{\mathcal
S}_n},{{\mathbb M}}) > t \bigr) \leq \frac{2^d}{\rho(t/2)
t^d}\exp \bigl( -n
\rho(t) t^d \bigr),
\end{equation}
where $\rho(t)$ is defined in~(\ref{eqrho-t}).
If, in addition, $t < \min\{ \rho/(2C_2), t_0\}$, then
%
\begin{equation}
\mathbb{P}\bigl(\mm{\mathbf{H}}({{\mathcal S}_n},{{\mathbb M}}) > t
\bigr) \leq \frac
{2^{d+1}}{t^d \rho
}\exp \biggl( -n \frac{\rho t^d}{2} \biggr).
\end{equation}
Hence if
$t_n(\alpha)< \min\{ \rho/(2C_2), t_0\}$ is the solution to the equation
%
\begin{equation}
\label{eqlambert} \frac{2^{d+1}}{t_n^d \rho}\exp \biggl( -n \frac{\rho t_n^d}{2} \biggr)=
\alpha,
\end{equation}
then
\[
\mathbb{P} \bigl(W_\infty(\widehat{\mathcal P},{\mathcal P}) >
t_n(\alpha) \bigr) \leq \mathbb{P} \bigl(\mm{\mathbf{H}}({{\mathcal
S}_n},{{\mathbb M}}) > t_n(\alpha) \bigr) \leq\alpha.
\]
\end{lemma}

\begin{rems*}
From the previous lemma,
it follows that, setting $t_n =\break  ( \frac{4}{\rho} \frac{\log n}{n}
 )^{1/d}$,
\[
\mathbb{P} \bigl(\mm{\mathbf{H}}({{\mathcal S}_n},{{\mathbb M}}) >
t_n \bigr) \leq\frac{2^{d-1}}{n \log n}, %
\]
for all $n$ large enough.
The right-hand side of~(\ref{eqlam1})
is known as the Lambert function
[\citet{lambert}].
Equation~(\ref{eqlambert})
does not admit a closed form solution, but
can be solved numerically.

To use the lemma, we need to
estimate $\rho$.
Let $P_n$ be the empirical measure induced by the sample~$\mathcal
{S}_n$, given
by
\[
P_n(A) = \frac{1}{n}\sum_{i=1}^n
I_A(X_i), %
\]
for any measurable Borel set $A \subset\mathbb{R}^D$.
Let $r_n$ be a positive small number and
consider the plug-in estimator of $\rho$,
%
\begin{equation}
\hat{\rho}_n = \min_i \frac{P_n(B(X_i,r_n/2))}{r_n^d}.
\end{equation}

Our next result shows that, under our assumptions and provided that
the sequence $r_n$ vanishes at an appropriate rate as
$n \to\infty$, $\hat{\rho}_n$
is a consistent estimator of $\rho$.
\end{rems*}

%
\begin{theorem}
\label{theoremahat}
Let
$r_n \asymp(\log n /n)^{1/(d+2)}$.
Then
\[
\hat{\rho}_n - \rho= O_P ( r_n ).
\]
\end{theorem}

\begin{rem*}
We have assumed that $d$ is known.
It is also possible to estimate $d$,
although we do not pursue that extension here.

We now need to use $\hat{\rho}_n$ to estimate $t_n(\alpha)$ as follows.
Assume that $n$ is even, and split the data randomly into two halves,
$\mathcal{S}_n = \mathcal{S}_{1,n} \sqcup\mathcal{S}_{2,n}$. Let
$\hat{\rho}_{1,n}$ be the plug-in estimator of $\rho$ computed from
$\mathcal{S}_{1,n}$, and define ${\hat t}_{1,n}$ to solve the equation
%
\begin{equation}
\label{eqsolvetn} \frac{2^{d+1}}{t^d \hat{\rho}_{1,n}} \exp \biggl( - \frac{n
t^d
\hat{\rho}_{1,n}}{2} \biggr) =
\alpha.
\end{equation}
\end{rem*}


%
\begin{theorem}
\label{theoremahat2}
Let $\widehat{\mathcal{P}}_2$ be the persistence diagram for the distance
function to~$\mathcal{S}_{2,n}$, then
%
\begin{eqnarray}
\mathbb{P}\bigl( W_{\infty}(\widehat{\mathcal{P}}_2,
\mathcal{P}) > \hat{t}_{1,n}\bigr) &\leq& \mathbb{P}\bigl(\mathbf{H}({
\mathcal S}_{2,n},{\mathbb M}) > \hat t_{1,n}\bigr)
\nonumber\\[-8pt]\\[-8pt]\nonumber
&\leq&
\alpha+ O \biggl(\frac{ \log n}{n} \biggr)^{1/(2+d)},\nonumber
\end{eqnarray}
where the probability $\mathbb{P}$ is with respect to both the joint
distribution of the entire sample and
the randomness induced by the sample splitting.
\end{theorem}

In practice, we have found that solving (\ref{eqsolvetn}) for
$\hat t_n$
without
splitting the data also works well although we do not have a formal proof.
Another way to define $\hat t_n$ which is simpler but more conservative,
is to define
%
\begin{equation}
\hat t_n = \biggl( \frac{2}{n\hat{\rho}_n} \log \biggl(
\frac{n}{\alpha
} \biggr) \biggr)^{1/d}.
\end{equation}
Then
$\hat t_n = u_n (1+ O(\hat{\rho}_n-\rho))$
where
$u_n =  ( \frac{\rho}{n\hat{\rho}_n} \log ( \frac
{n}{\alpha} )
 )^{1/d}$, and so
\begin{eqnarray*}
\mathbb{P}\bigl(\mathbf{H}({\mathcal S}_n,{\mathbb M}) > \hat t_n\bigr) &=& \mathbb{P}\bigl(\mathbf{H}({\mathcal S}_n,{
\mathbb M}) > u_n\bigr) + O \biggl(\frac{\log n}{n}
\biggr)^{1/(2+d)}
\\
&\leq&\alpha+ O \biggl(\frac{\log n}{n}
\biggr)^{1/(2+d)}.
\end{eqnarray*}

\subsection{Method \textup{III}: The method of shells}

The dependence of the previous method
on the parameter $\rho$
makes the method very fragile.
If the density is low in even a small region, then
the method above is a disaster.
Here we develop a sharper bound
based on shells of the form
$\{x\dvtx  \gamma_j < \rho(x, \downarrow0) < \gamma_{j+1}\}$
where we recall from~(\ref{eqrho-t}) and~(\ref{eqrho}) that
\[
\rho(x, \downarrow0) = \lim_{t \to0} \rho(x,t) = \lim
_{t \to0} \frac{P(B(x,t/2))}{t^d}.
\]

Let
$G(v) = P( \rho(X, \downarrow0) \leq v)$, and let $g(v) = G'(v)$.

%
\begin{theorem}
\label{theoremshells}
Suppose that $g$ is bounded and has a uniformly bounded, continuous derivative.
Then
for any $t\leq\rho/(2C_1)$,
%
\begin{eqnarray}
\quad\mathbb{P}\bigl( W_\infty(\widehat{\mathcal P},{\mathcal P}) > t\bigr)
&\leq& P\bigl( \mm{\mathbf{H}}({{\mathcal S}_n},{{\mathbb M}}) > t\bigr)
\nonumber\\[-8pt]\\[-8pt]\nonumber
&\leq& \frac{2^{d+1}}{t^d}\int_{\rho}^{\infty}
\frac{g(v)}{v} e^{-n v
t^d/2} \,dv.
\end{eqnarray}
\end{theorem}

Let $K$ be a smooth, symmetric kernel
satisfying the
conditions in \citet{gine2002rates}
(which includes all the commonly used kernels),
and let
%
\begin{equation}
\hat g(v) = \frac{1}{n}\sum_{i=1}^n
\frac{1}{b} K \biggl( \frac{v - V_i}{b} \biggr),
\end{equation}
where $b>0$, $V_i = \hat\rho(X_i,r_n)$, and
\[
\hat\rho(x,r_n) = \frac{P_n(B(x,r_n/2))}{r_n^d}. %
\]

%
\begin{theorem}
\label{theoremghat}
Let
$r_n =  (\frac{ \log n}{n} )^{1/(d+2)}$.
\begin{longlist}[(1)]
\item[(1)] We have that
\[
\sup_v \bigl\llvert \hat g(v) - g(v)\bigr\rrvert =
O_P \biggl( b^2 + \sqrt{\frac{\log n}{nb}}+
\frac{r_n}{b^2} \biggr). %
\]
Hence if we choose
$b\equiv b_n \asymp r_n^{1/4}$,
then
\[
\sup_v \bigl\llvert \hat g(v) - g(v)\bigr\rrvert =
O_P \biggl(\frac{\log n}{n} \biggr)^{1/(2(d+2))}. %
\]

\item[(2)]
Suppose that $n$ is even and that $\rho>0$.
Assume that $g$ is strictly positive over its support
$[\rho,B]$.
Randomly split the data into two halves: $\mathcal{S}_n =
(\mathcal{S}_{1,n},\mathcal{S}_{2,n})$.
Let $\hat{g}_{1,n}$ and $\hat{\rho}_{1,n}$ be estimators of
$g$ and $\rho$,
respectively, computed from the first half of the data, and
define $\hat{t}_{1,n}$ to be the solution of the equation
%
\begin{equation}
\frac{2^{d+1}}{\hat{t}{}^d_{1,n}}\int_{\hat{\rho}_n}^\infty
\frac{\hat g(v)}{v} e^{-n v \hat{t}{}^d/2} \,dv = \alpha.
\end{equation}
Then
%
\begin{equation}
\mathbb{P}\bigl( W_{\infty}(\widehat{\mathcal{P}}_2,
\mathcal{P}) > \hat{t}_{1,n}\bigr) \leq \mathbb{P}\bigl(\mathbf{H}({
\mathcal S}_{2,n},{\mathbb M}) > \hat t_{1,n}\bigr) \leq
\alpha+ O(r_n),
\end{equation}
where $\widehat{\mathcal{P}}_2$ is the persistence diagram associated to
$\mathcal{S}_{2,n}$ and
the probability $\mathbb{P}$ is with respect to both the joint
distribution of
the entire sample and
the randomness induced by the sample splitting.
\end{longlist}
\end{theorem}

\subsection{Method \textup{IV}: Density estimation}\label{sectiondensity}

\begin{figure}[b]

\includegraphics{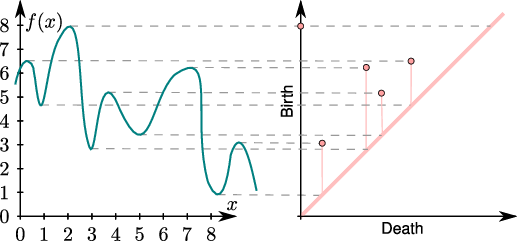}

\caption{We plot the persistence diagram corresponding to
the upper level
set filtration of the function $f(x)$. For consistency,
we swap the birth and death axes so all persistence points
appear above the line $y=x$. The point born first does not
die, but for convenience we mark it's death at $f(x)=0$.
This is analogous to the point marked with $\infty$ in
our previous diagrams.}
\label{figupperLevelPers}
\end{figure}

In this section, we take a completely different approach.
We use the data to construct a smooth density estimator, and then
we find the persistence diagram
defined by a filtration of the upper level sets
of the density estimator; see~Figure~\ref{figupperLevelPers}.
A different approach to smoothing based on diffusion distances is
discussed in \citet{bendich2011improving}.

%

Again, let $X_1,\ldots, X_n$ be a sample from $P$.
Define
%
\begin{equation}
p_h(x) = \int_{{\mathbb M}} \frac{1}{h^D} K \biggl(
\frac{\llVert  x-u\rrVert  _2}{h} \biggr) \,dP(u).
\end{equation}
Then\vspace*{1pt} $p_h$ is the density of the probability measure
$P_h$ which is the convolution
$P_h = P\star\mathbb{K}_h$
where
$\mathbb{K}_h(A) = h^{-D} \mathbb{K}(h^{-1} A)$ and
$\mathbb{K}(A) = \int_A K(t) \,dt$.
That is, $P_h$
is a smoothed version of $P$.
Our target of inference in this section is
${\mathcal P}_h$, the persistence diagram
of the upper level sets of $p_h$.
The standard estimator for $p_h$ is the kernel density estimator
%
\begin{equation}
\hat p_h(x) = \frac{1}{n}\sum
_{i=1}^n \frac{1}{h^D} K \biggl(
\frac{\llVert  x-X_i\rrVert  _2}{h} \biggr).
\end{equation}
It is easy to see that
$\mathbb{E}(\hat p_h(x)) = p_h(x)$.
Let us now explain why
${\mathcal P}_h$ is of interest.

First, the upper level sets
of a density are of
intrinsic interest in statistics in machine learning.
The connected components of the upper level sets
are often used for clustering.
The homology of these upper level sets provides
further structural information about the density.

Second, under appropriate conditions,
the upper level sets of $p_h$
may carry topological information about a set
of interest ${\mathbb M}$.
To see this, first suppose that
${\mathbb M}$ is a smooth, compact $d$-manifold,
and suppose that $P$ is supported on ${\mathbb M}$.
Let $p$ be the density of $P$
with respect to Hausdorff measure on ${\mathbb M}$.
In the special case where $P$ is the uniform distribution,
every upper level set $\{p>t\}$ of $p$
is identical to ${\mathbb M}$ for $t>0$ and small enough.

Thus if ${\mathcal P}$ is the persistence diagram defined
by the upper level sets
$\{x\dvtx  p(x) > t\}$
of $p$,
and ${\mathcal Q}$ is the persistence diagram of the distance function
$d_{{\mathbb M}}$, then the points of ${\mathcal Q}$ are in one-to-one
correspondence with the generators of $H({\mathbb M})$, and
the points with
higher persistence in ${\mathcal P}$
are also in 1--1 correspondence with the generators of $H({\mathbb M})$.
For example, suppose that
${\mathbb M}$ is a circle in the plane with radius $\tau$.
Then ${\mathcal Q}$ has two points: one at
$(0,\infty)$ representing a single connected component, and one at
$(0,\tau)$ representing the single cycle.
${\mathcal P}$ also has two points: both at
$(0,1/2\pi\tau)$
where $1/2\pi\tau$ is simply the maximum of the density over the circle.
In sum, these two persistence diagrams
contain the same information; furthermore,
$\{x\dvtx  p(x) > t\} \cong{\mathbb M}$ for all
$0< t < 1/2\pi\tau$.

If $P$ is not uniform but has a smooth density $p$, bounded away from 0,
then there is an interval $[a,A]$ such that
$\{x\dvtx  p(x) > t\} \cong{\mathbb M}$
(i.e., is homotopic)
for
$a \leq t \leq A$.
Of course, one can create examples where no level sets
are equal to ${\mathbb M}$, but it seems unlikely that any
method can recover the homology of ${\mathbb M}$ in those cases.

Next, suppose there is noise; that is,
we observe
$Y_1,\ldots, Y_n$, where $
Y_i = X_i + \sigma\varepsilon_i
$
and $\varepsilon_1,\ldots,\varepsilon_n \sim\Phi$.
We assume that
$X_1,\ldots, X_n \sim Q$
where $Q$ is supported on ${\mathbb M}$.
Note that
$X_1,\ldots, X_n$ are unobserved.
Here, $\Phi$ is the noise distribution and $\sigma$
is the noise level.
The distribution $P$ of $Y_i$ has density
%
\begin{equation}
p(y) = \int_{{\mathbb M}} \phi_\sigma(y-u) \,dQ(u),
\end{equation}
where
$\phi$ is the density of $\varepsilon_i$ and
$\phi_\sigma(z) = \sigma^{-D}\phi(y/\sigma)$.
In this case,
no level set
$L_t = \{y\dvtx  p(y) > t\}$
will equal ${\mathbb M}$.
But as long as $\phi$ is smooth and $\sigma$ is small,
there will be a range of values
$a\leq t \leq A$ such that
$L_t\cong{\mathbb M}$.

The estimator
$\hat p_h(x)$
is consistent for $p$
if $p$ is continuous,
as long as we let the bandwidth $h=h_n$
change with $n$
in such a way that
$h_n \to0$
and $n h_n \to\infty$
as $n\to\infty$.
However, for exploring topology,
we need not let $h$ tend to 0.
Indeed, more precise topological inference is obtained by using
a bandwidth $h>0$.
Keeping $h$ positive smooths the density,
but the level sets can still retain the correct topological information.

We would also like to point out that
the quantities
${\mathcal P}_h$ and
$\widehat{\mathcal P}_h$ are
more robust and
much better
behaved statistically
than the \v{C}ech complex of the raw data.
In the language of computational topology,
${\mathcal P}_h$ can be considered a topological simplification of
${\mathcal P}$.
${\mathcal P}_h$ may omit subtle details that are present in ${\mathcal P}$
but is much more stable.
For these reasons,
we now focus on estimating
${\mathcal P}_h$.

Recall that, from the stability theorem,
%
\begin{equation}
W_\infty(\widehat{\mathcal P}_h,{\mathcal
P}_h) \leq\llVert \hat p_h - p_h
\rrVert _\infty.
\end{equation}
Hence it suffices to find $c_n$ such that
%
\begin{equation}
\limsup_{n\to\infty}\mathbb{P}\bigl(\llVert \hat p_h
- p_h\rrVert _\infty> c_n\bigr) \leq\alpha.
\end{equation}

\subsubsection*{Finite sample band}
Suppose that the support of $P$ is
contained in ${\mathcal X} = [-C,C]^D$.
Let $p$ be the density of $P$.
Let $K$ be a kernel with the same assumptions as above, and choose a bandwidth
$h$. Let
\[
\hat p_h(x) = \frac{1}{n} \sum
_{i=1}^n \frac{1}{h^D} K \biggl(
\frac{\llVert  x-X_i\rrVert  _2}{h} \biggr) %
\]
be the kernel density estimator, and let
\[
p_h(x) = \frac{1}{h^D} \int_{\mathcal X} K
\biggl(\frac
{\llVert  x-u\rrVert  _2}{h} \biggr) \,dP(u) %
\]
be the mean of $\hat p_h$.

%
\begin{lemma}
\label{lemmaHoeff}
Assume that $\sup_x K(x) = K(0)$ and
that $K$ is $L$-Lipschitz, that is,
$\llvert K(x)- K(y)\rrvert \leq L \llVert  x-y\rrVert  _2$.
Then
%
\begin{equation}
\mathbb{P} \bigl( \llVert \hat p_h - p_h\rrVert
_{\infty} > \delta \bigr) \leq 2 \biggl( \frac{4 C L \sqrt{D}}{\delta h^{D+1}}
\biggr)^D \exp \biggl( - \frac{ n \delta^2 h^{2D}}{2 K^2(0) } \biggr).
\end{equation}
\end{lemma}

\begin{rem*}
The proof of the above lemma uses Hoeffding's inequality.
A~sharper result can be obtained by using Bernstein's inequality;
however, this introduces extra constants that need to be estimated.

We\vspace*{2pt} can use the above lemma to approximate the persistence diagram
for $p_h$, denoted by $\mathcal{P}_h$, with the diagram
for $\hat{p}_h$, denoted by $\widehat{\mathcal{P}}_h$:
\end{rem*}

%
\begin{corollary}
Let $\delta_n$ solve
%
\begin{equation}
2 \biggl( \frac{4 C L \sqrt{D}}{\delta_n h^{D+1}} \biggr)^D \exp \biggl( -
\frac{ n \delta_n^2 h^{2D}}{2 K^2(0) } \biggr) = \alpha.
\end{equation}
Then
%
\begin{equation}
\sup_{P\in{\mathcal Q}} \mathbb{P} \bigl( W_\infty(\widehat {
\mathcal P}_h,{\mathcal P}_h) > \delta_n
\bigr) \leq \sup_{P\in{\mathcal Q}} \mathbb{P} \bigl( \llVert \hat p_h - p_h\rrVert _{\infty} > \delta_n
\bigr) \leq \alpha,
\end{equation}
where ${\mathcal Q}$ is the set of all probability measures supported
on ${\mathcal X}$.
\end{corollary}

Now we consider a different finite sample band.
Computationally, the persistent homology of the upper level sets of
$\hat p_h$
is actually based on a piecewise linear approximation to $\hat p_h$.
We choose a finite grid $G\subset\mathbb{R}^D$
and form a triangulation over the grid.
Define $\hat p_h^\dagger$ as follows.
For $x\in G$, let
$\hat p_h^\dagger(x) = \hat p_h(x)$.
For $x\notin G$, define
$\hat p_h^\dagger(x)$ by linear interpolation over the triangulation.
Let
$p_h^\dagger(x)= \mathbb{E}(\hat p_h^\dagger(x))$.
The real object of interest is the persistence diagram
${\mathcal P}_h^\dagger$ of the upper level set filtration of
$p_h^\dagger(x)$. We
approximate this diagram with the persistence diagram $\widehat
{\mathcal P}_h^\dagger$
of the upper level set filtration of $\hat{p}_h^\dagger(x)$.
As before,
\[
W_\infty\bigl(\widehat{\mathcal P}_h^\dagger,{
\mathcal P}_h^\dagger\bigr) \leq \bigl\llVert \hat p_h^\dagger- p_h^\dagger\bigr\rrVert
_\infty. %
\]
But due to the piecewise linear nature of these functions,
we have that
\[
\bigl\llVert \hat p_h^\dagger- p_h^\dagger
\bigr\rrVert _\infty\leq\max_{x\in G} \bigl\llvert
\hat p_h^\dagger(x) - p_h^\dagger(x)
\bigr\rrvert. %
\]

%
\begin{lemma}
Let $N = \llvert G\rrvert $ be the size of the grid.
Then
%
\begin{equation}
\mathbb{P} \bigl(\bigl\llVert \hat p_h^\dagger-
p_h^\dagger\bigr\rrVert _{\infty} > \delta \bigr) \leq
2 N \exp \biggl( - \frac{ 2n \delta^2 h^{2D}}{K^2(0) } \biggr).
\end{equation}
Hence, if
%
\begin{equation}
\delta_n = \biggl(\frac{K(0)}{h} \biggr)^D \sqrt{
\frac{1}{2n} \log \biggl( \frac{2N}{\alpha} \biggr),}
\end{equation}
then
%
\begin{equation}
\mathbb{P} \bigl(\bigl\llVert \hat p_h^\dagger-
p_h^\dagger\bigr\rrVert _{\infty} >
\delta_n \bigr) \leq\alpha.
\end{equation}
\end{lemma}

This band can be substantially tighter as long we do not use
a grid that is too fine.
In a sense, we are rewarded for acknowledging that
our topological inferences take place
at some finite resolution.

\subsubsection*{Asymptotic confidence band}
A tighter---albeit only asymptotic---bound
can be obtained using large sample theory.
The simplest approach is the bootstrap.

Let $X_1^*,\ldots, X_n^*$ be a sample
from the empirical distribution $P_n$, and
let $\hat p^*_h$ denote the density estimator
constructed from
$X_1^*,\ldots, X_n^*$.
Define the random measure
%
\begin{equation}
J_n(t) = \mathbb{P} \bigl( \sqrt{n h^D}\bigl\llVert
\hat p_h^* - \hat p_h\bigr\rrVert
_\infty> t | X_1,\ldots, X_n \bigr)
\end{equation}
and the bootstrap quantile
$Z_\alpha= \inf\{ t\dvtx  J_n(t) \leq\alpha\}$.

%
\begin{theorem}
As $n\to\infty$,
\[
\mathbb{P} \biggl( W_\infty(\widehat{\mathcal P}_h,{
\mathcal P}_h) > \frac{Z_\alpha}{\sqrt{nh^D}} \biggr) \leq \mathbb{P} \bigl(
\sqrt{nh^D}\llVert \hat p_h - p_h\rrVert
_\infty> Z_\alpha \bigr) = \alpha+ O \biggl(\sqrt{
\frac{1}{n}} \biggr). %
\]
\end{theorem}

The proof follows from standard results; see, for example,
\citet{chazal2013bootstrap}.
As usual, we approximate $Z_\alpha$ by Monte Carlo.
Let
$T = \sqrt{n h^D}\llVert  \hat p_h - \hat p_h^*\rrVert  _\infty$
be from a bootstrap sample.
Repeat bootstrap $B$ times
yielding values
$T_1,\ldots, T_B$.
Let
\[
Z_\alpha= \inf \Biggl\{ z\dvtx  \frac{1}{B}\sum
_{j=1}^B I( T_j > z) \leq \alpha \Biggr
\}. %
\]
We can ignore the error due to the fact that $B$ is finite
since this error can be made as small as we~like.

\begin{rem*}
We have emphasized fixed $h$ asymptotics
since, for topological inference, it is not necessary to let $h\to0$
as $n\to\infty$.
However, it is possible to let $h\to0$ if one wants.
Suppose $h\equiv h_n$ and $h\to0$ as $n\to\infty$.
We require that
$n h^D/\log n \to\infty$ as $n\to\infty$.
As before, let $Z_\alpha$ be the bootstrap quantile.
It follows from Theorem~3.4 of
\citet{neumann1998strong},
that
%
\begin{eqnarray}
\mathbb{P} \biggl( W_\infty(\widehat{\mathcal P}_h,{
\mathcal P}_h) > \frac{Z_\alpha}{\sqrt{nh^D}} \biggr) &\leq& \mathbb{P} \bigl(
\sqrt{n h^D}\llVert \hat p_h - p_h
\rrVert _\infty> Z_\alpha \bigr)
\nonumber\\[-8pt]\\[-8pt]\nonumber
&=& \alpha + \biggl(
\frac{\log n}{n h^D} \biggr)^{(4+D)/(2(2+D))}.
\end{eqnarray}
\end{rem*}

\begin{ou*}
Now
%
%
\begin{figure}

\includegraphics{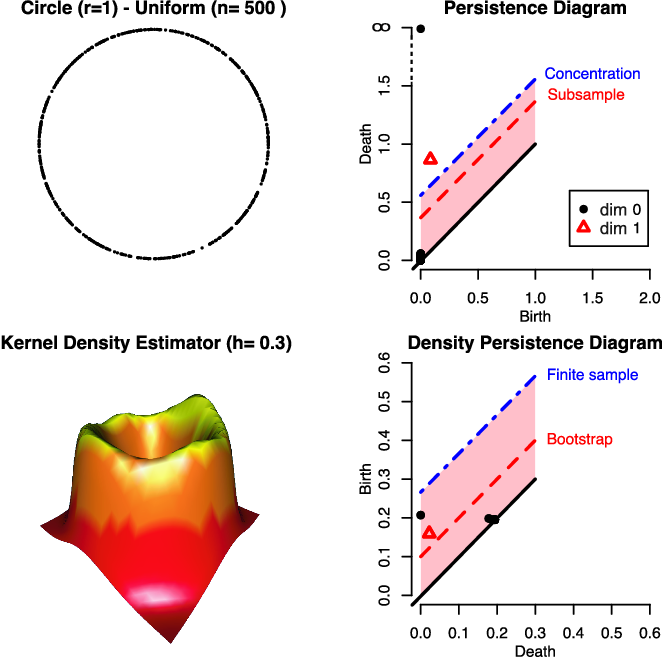}

\caption{Uniform distribution over the unit Circle. (Top left) sample ${\mathcal
S}_n$. (Top right)~corresponding persistence diagram.
The black
circles indicate the life span of connected components, and the red
triangles indicate the life span of 1-dimensional holes. (Bottom left) kernel
density estimator. (Bottom right) density persistence diagram. For more details
see Example \protect\ref{exCircleUniform}.}\label{figCircleUniform}
\end{figure}
we explain why the density-based method is very insensitive to outliers.
Let
$P= \pi U + (1-\pi) Q$
where $Q$ is supported on ${\mathbb M}$,
$\pi>0$ is a small positive constant and
$U$ is a smooth distribution supported on
$\mathbb{R}^D$.
Apart from a rescaling,
the bottleneck distance between
${\mathcal P}_P$ and
${\mathcal P}_Q$ is at most $\pi$.
The kernel estimator is still a consistent estimator of $p$, and
hence the persistence diagram is barely affected by outliers.
In fact, in the examples in
\citet{bendich2011improving},
there are only a few outliers which formally corresponds to
letting $\pi= \pi_n \to0$ as $n\to\infty$.
In this case, the density method is very robust.
We show this in more detail in the experiments section.
\end{ou*}

\section{Experiments}\label{sectionexperiments}

As is common in the literature
on computational topology, we focus
on a few simple, synthetic examples.
For each of them we compute the Rips persistence diagram and the
density persistence diagram introduced in Section~\ref{sectiondensity}. We use a Gaussian kernel with bandwidth $h=0.3$.
This will serve to illustrate the different methods for the
construction of
confidence bands for the persistence diagrams.

%
\begin{example}
\label{exCircleUniform}
Figure~\ref{figCircleUniform} shows
%
%
\begin{figure}

\includegraphics{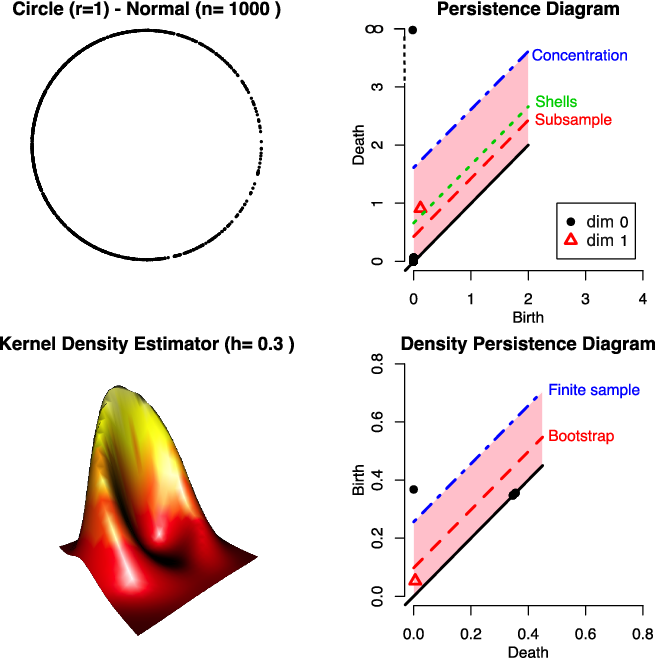}

\caption{Truncated Normal distribution over the unit Circle. (Top left) sample
${\mathcal S}_n$. (Top right) corresponding persistence diagram. The
black circles
indicate the life span of connected components, and the red triangles indicate
the life span of 1-dimensional holes. (Bottom left) kernel density estimator.
(Bottom right) density persistence diagram. For more details see
Example \protect\ref{exCircleNormal}.}\label{figCircleNormal}
\end{figure}
the methods described in the
previous sections applied to a sample
from the uniform distribution over the unit circle ($n=500$). In the top
right plot the
different 95\% confidence bands for the persistence diagram are
computed using methods~I~(subsampling) and~II (concentration of measure). Note that the
uniform distribution does not satisfy the assumptions for the method
of shells.
The subsampling method
and the concentration method both correctly
show one significant connected component
and one significant loop.
In the bottom right plot the finite sample density estimation method
and the
bootstrap method are applied to the density persistence diagram. The first
method does not have sufficient power to detect the topological features.
However, the bootstrap
density estimation method does find that one connected component and
one loop
are significant.\looseness=-1
\end{example}

%
\begin{example}
\label{exCircleNormal}
Figure~\ref{figCircleNormal} shows
%
%
\begin{figure}

\includegraphics{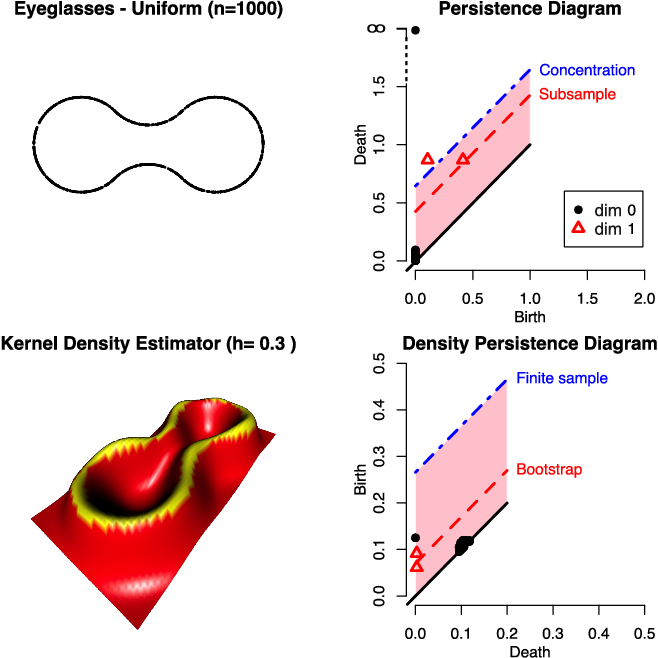}\vspace*{-1pt}

\caption{Uniform distribution over the eyeglasses curve. (Top left)
sample ${\mathcal
S}_n$. (Top right) corresponding persistence diagram. The black circles indicate
the life span of connected components and the red triangles indicate
the life
span of 1-dimensional holes. Bottom left: kernel density estimator. (Bottom
right) density persistence diagram. For more details see Example
\protect\ref{exGlassesUniform}.}\vspace*{-5pt}\label{figGlassesUniform}
\end{figure}
the methods described in the
previous sections applied to a sample from the truncated Normal distribution
over the unit circle ($n=1000$). The top left plot shows the sample, and the
bottom left plot shows the kernel density estimator constructed using
a Gaussian kernel with bandwidth $h=0.3$.
The plots on the left show the different methods for the construction
of 95\% confidence bands around the diagonal of the persistence diagrams.
This case is challenging because there is a portion of the circle that is
sparsely sampled.
The concentration method fails to detect the loop, as shown in the top right
plot.
However, the method of shells and the subsampling method both
declare the loop to be significant.
The bottom right plot shows that both the finite sample method and the bootstrap
method
fail to detect the loop.
\end{example}

%
\begin{example}
\label{exGlassesUniform}
Figure~\ref{figGlassesUniform} shows
%
%
\begin{figure}

\includegraphics{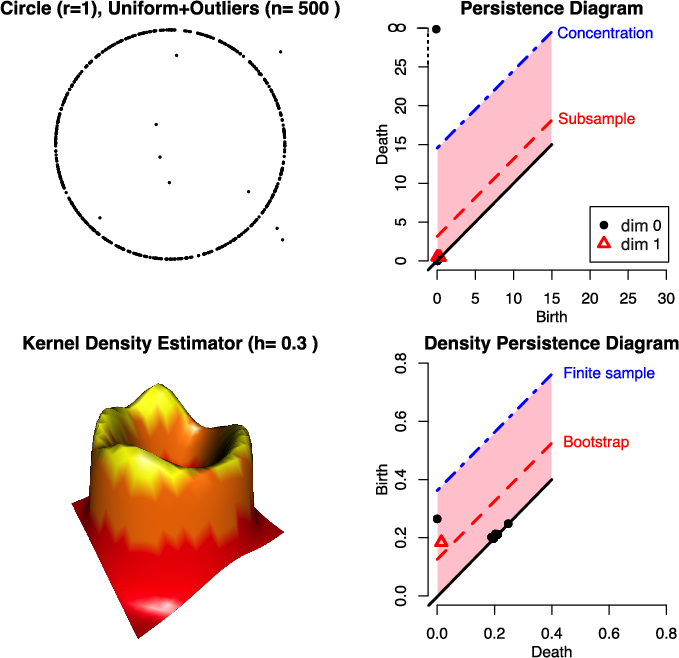}

\caption{Uniform distribution over the unit Circle with some outliers.
See Example \protect\ref{exOutliers}.}
\label{figCircleUniformOut}
\end{figure}
the methods described in the
previous sections applied to a sample of size
$n=1000$ from the uniform distribution over the eyeglasses, a figure
similar to the Cassini Curve obtained by gluing together two unit circles.
Note that the uniform distribution does not satisfy the assumptions for the
method
of shells.
Each method provides a 95\% confidence band around the diagonal for
the persistence diagrams.
The top right plot shows that the subsample method declares both the
loops to
be significant, while the concentration method detects only one of
them, as
the bootstrap method for the density persistence diagram, shown in the
bottom right plot.\vspace*{-2pt}
\end{example}

\begin{figure}

\includegraphics{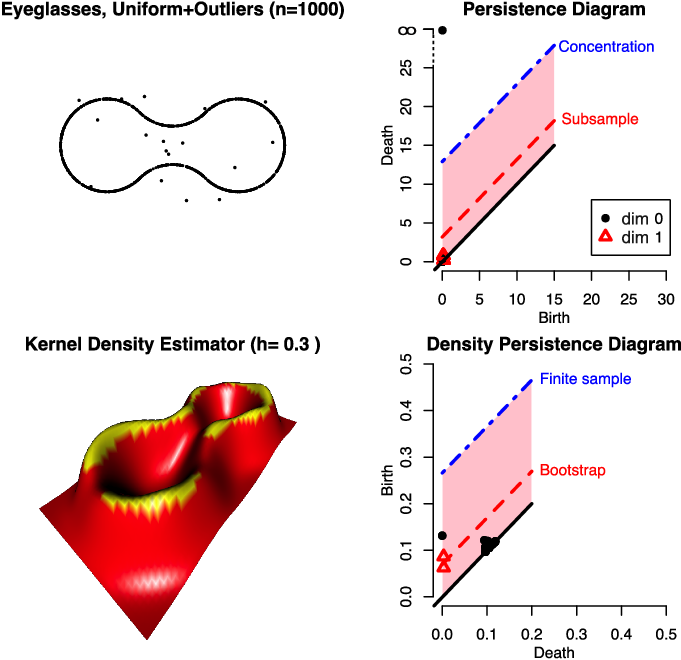}

\caption{Uniform distribution over the eyeglasses curve with some
outliers. See Example~\protect\ref{exOutliers}.}\vspace*{-5pt}\label{figGlassesUniformOut}
\end{figure}

%
\begin{example}
\label{exOutliers}
In this example, we replicate Examples~\ref{exCircleUniform} and
\ref{exGlassesUniform}, adding some outliers to the uniform
distributions over
the unit circle and the\vadjust{\goodbreak} eyeglasses.
Figures~\ref{figCircleUniformOut} and~\ref{figGlassesUniformOut}
show the
persistence
diagrams with different methods for the construction of 95\% confidence bands.
Much of the literature on computational topology focuses
on methods that use the distance function to the data.
As we see here,
and as discussed in
\citet{bendich2011improving},
such methods are quite fragile. A few outliers are sufficient
to drastically change the persistence diagram and force
the concentration method and the subsample method
to declare the topological features to be not significant.
On the other hand the density-based methods are very insensitive
to the presence of outliers, as shown in the bottom right plots
of the two figures.
\end{example}

\section{Proofs}\label{sectionproofs}
In this section, we provide proofs of the theorems and lemmas found in
Section~\ref{sectionmodel}.

%

Recall that the $\delta$-covering number $N$ of a manifold ${\mathbb
M}$ is the smallest
number of Euclidean balls of radius $\delta$ required to cover the
set. The
$\delta$-packing number $N'$ is
the maximum number of sets of the form $B(x,\delta) \cap{\mathbb M}$,
where $x \in{\mathbb M}$,
that may be packed into ${\mathbb M}$ without overlap.
First, we prove the following lemma.

%
\begin{lemma}
\label{lemmapacking}
Let $r_n=o(1)$ and let $\mathcal{D}_n = \{d_1, \ldots, d_{N'} \}$
be the set of centers of $\{B(d_i,r_n) \cap{\mathbb M}\dvtx  i=1,\ldots, N'\}$, a
$r_n$-packing set for~$\mathbb{M}$ (this set is nonempty and of cardinality
$N'$ increasing in $n$). For large $n$, the size of the $r_n$-packing satisfies
%
\begin{equation}
\label{eqNprime} \frac{\operatorname{vol}({\mathbb M})}{A_2(d) r_n^d} \leq N'\leq \frac{A_1(d)} {r_n^d},
\end{equation}
for some constants $A_1(d)$ and $A_2(d)$ depending on the reach and
dimension $d$ of~${\mathbb M}$.
\end{lemma}
\begin{pf}
From the proof of Lemma~\ref{lemmaconcentration} we know that
%
\begin{equation}
N'\leq\frac{A_1(d)} {r_n^d}, \label{equpperN}
\end{equation}
where $A_1(d)$ is a constant depending on the dimension of the manifold.
A similar lower bound for $N'$ can be obtained. Let $N$ be the size of a
$2r_n$-covering set for $\mathbb{M}$, formed by Euclidean balls $B(c_i,2r_n)$
with centers $\mathcal{C}= \{c_1, \ldots, c_{N} \}$. By Lemma~5.2 in
\citet{niyogi2008finding} and a simple volume argument we have
%
\begin{equation}
N' \geq N \geq\frac{\operatorname{vol}({\mathbb M})}{\max_{i=1,\ldots,N}
\operatorname{vol} (B(c_i,2r_n) \cap{\mathbb M} )}. \label{eqlowerN1}
\end{equation}
For large $n$, by Corollary 1.3 in \citet{Fred2013Note},
%
\begin{equation}
\max_{i=1,\ldots,N} \operatorname{vol} \bigl(B(c_i,2r_n)
\cap{\mathbb M} \bigr) \leq A_2(d) r_n^d,
\label{eqlowerN2}
\end{equation}
where $A_2(d)$ is a constant depending on the dimension of the manifold.
Combining~(\ref{equpperN}), (\ref{eqlowerN1}) and (\ref{eqlowerN2})
we obtain (\ref{eqNprime}).
\end{pf}

\begin{pf*}{Proof of Theorem \protect\ref{theoremsubsample}}
We begin the proof by showing that there exists an event of probability
approaching one (as $n \rightarrow\infty$) such
that, over this event, for any subsample
${\mathcal S}_{b,n}$ of the data $\mathcal{S}_n$ of size $b$,
%
\begin{equation}
\label{eqSnSbn} \mm{\mathbf{H}}({{\mathcal S}_n},{{\mathbb M}})\leq
\mm{\mathbf {H}}({{\mathcal S}_{b,n}},{{\mathcal S}_n}).
\end{equation}
Toward this end,
let~$t_n =  ( \frac{4}{\rho} \frac{\log n}{n}  )^{1/d}$, and
define the event $\mathcal{A}_n = \{ \mm{\mathbf{H}}({\mathcal
{S}_n},{{\mathbb M}}) < t_n\}$.
Then, by the remark following\vspace*{2pt} Lemma~\ref{lemmaconcentration},
$\mathbb{P}(\mathcal{A}_n^c) \leq\frac{2^{d-1}}{n \log n}$, for all
$n$ large
enough.
Next, let $\mathcal{D}_n = \{d_1, \ldots, d_{N'} \}$
be the set of centers of $\{B(d_i,2t_n) \cap{\mathbb M}\dvtx  i=1,\ldots, N'\}$, a
$2t_n$-packing set for~$\mathbb{M}$. By Lemma~\ref{lemmapacking},
the size of
$\mathcal{D}_n$ is of order $t_n^{-d} = \Theta (
\frac{n}{\log n}  )$.

We can now show (\ref{eqSnSbn}). Suppose $\mathcal{A}_n$ holds and,
arguing by
contradiction,
also that
$\mm{\mathbf{H}}({{\mathcal S}_{b,n}},{{\mathcal S}_n}) < \mm
{\mathbf{H}}({{\mathcal S}_n},{{\mathbb M}})$. Then
%
\begin{eqnarray}\label{eq2tn}
\mm{\mathbf{H}}({{\mathcal S}_{b,n}},{
\mathcal{D}_n}) &\leq&\mm {\mathbf{H}}({{\mathcal S}_{b,n}},{{
\mathcal S}_n}) + \mm{\mathbf{H}}({{\mathcal S}_n},{
\mathcal{D}_n}) < \mm{\mathbf{H}}({{\mathcal S}_n},{{
\mathbb M}}) + \mm{\mathbf {H}}({{\mathcal S}_n},{\mathcal{D}}) \nonumber
\\
&\leq& 2 \mm{\mathbf {H}}({{\mathcal S}_n},{{\mathbb M}})
\\
&\leq& 2t_n.\nonumber
\end{eqnarray}
Because of our assumption on $b$ and of (\ref{eqNprime}), $\frac{b}{N'}
\rightarrow0$ as
$n \rightarrow\infty$, which implies that a $(1-o(1))$ fraction of
the balls $\{ B(d_j,2t_n), j=1,\ldots,N' \}$ contains no points from
${\mathcal
S}_{b,n}$.
So
$\mm{\mathbf{H}}({{\mathcal S}_{b,n}},{\mathcal{D}}) > 2t_n$,
which, in light of~(\ref{eq2tn}), yields a contradiction. Thus, on
$\mathcal{A}_n$, (\ref{eqSnSbn}) holds,
for any subsample $\mathcal{S}_{b,n}$ of size $b$,
as claimed.

Next, let $ \{ \mathcal{S}_{b,n}^j, j=1,\ldots,N \}$ be an
enumeration of all
possible
subsamples of $\mathcal{S}_n$ of size $b$, where $N = {n \choose b}$, and
define
\[
\widetilde{L}_b(t) =\frac{1}{N}\sum
_{j=1}^N I\bigl(\mm{\mathbf{H}}\bigl({
\mathcal{S}^j_{b,n}},{{\mathbb M}}\bigr) > t\bigr).
\]

Using (\ref{eqSnSbn}) we obtain that, on the event ${\mathcal A}_n$,
\[
\mm{\mathbf{H}}\bigl({{\mathcal S}^j_b},{{\mathbb M}}
\bigr) \leq\mm{\mathbf {H}}\bigl({{\mathcal S}^j_b},{{
\mathcal S}_n}\bigr) + \mm{\mathbf {H}}({{\mathcal
S}_n},{{\mathbb M}}) \leq2\mm{\mathbf {H}}\bigl({{\mathcal
S}^j_b},{{\mathcal S}_n}\bigr), %
\]
and therefore,
$\mm{\mathbf{H}}({{\mathcal S}^j_b},{{\mathcal S}_n})\geq\mm
{\mathbf{H}}({{\mathcal S}^j_b},{{\mathbb M}})/2$,
simultaneously over all $j =1,\ldots,N$.
Thus, on that event
%
\begin{equation}
\label{eqLb2} \widetilde{L}_b(t) \leq L_b(t/2)\qquad
\mbox{for all } t>0,
\end{equation}
where $L_b$ is defined in (\ref{eqLb1}).
Thus, letting $\mathbf{1}_{\mathcal{A}_n} =
\mathbf{1}_{\mathcal{A}_n}(\mathcal{S}_n)$ be the indicator function
of the
event $\mathcal{A}_n$, we obtain the bound
\[
\widetilde{L}_b(c_b) = \widetilde{L}_b(c_b)
( \mathbf{1}_{\mathcal
{A}_n} + \mathbf{1}_{\mathcal{A}^c_n}) \leq
L_b(c_b/2) + \mathbf{1}_{\mathcal{A}^c_n} \leq\alpha+
\mathbf{1}_{\mathcal{A}^c_n},
\]
where the first inequality is due to (\ref{eqLb2}) and the second
inequality to the fact that $L_b(c_b/2) = \alpha$, by definition of $c_b$.
Taking expectations, we obtain that
%
\begin{equation}
\label{eqLbcb} \mathbb{E} \bigl( \widetilde{L}_b(c_b)
\bigr) \leq\alpha+ \mathbb{P}\bigl( \mathcal{A}^c_n\bigr)
= \alpha+ O \biggl( \frac{1}{n \log n} \biggr).
\end{equation}

Next, define
\[
J_b(t) = \mathbb{P}\bigl(\mm{\mathbf{H}}({\mathcal{S}_b},{{
\mathbb M}}) > t\bigr), \qquad t > 0, %
\]
where we recall that $\mathcal{S}_b$ is an i.i.d. sample of size $b$.
Then Lemma A.2 in \citet{RomanoShaikh12} yields that, for any
$\epsilon> 0$,
%
\begin{equation}
\label{eqromano} \mathbb{P} \Bigl( \sup_{t > 0} \bigl\llvert
\widetilde{L}_b(t) - J_b(t)\bigr\rrvert > \epsilon
\Bigr) \leq\frac{1}{\epsilon} \sqrt{ \frac{2 \pi}{k_n } },
\end{equation}
where $k_n = \lfloor\frac{n}{b} \rfloor$.
Let $\mathcal{B}_{n}$ be the event that
\[
\sup_{t > 0} \bigl\llvert \widetilde{L}_b(t) -
J_b(t)\bigr\rrvert \leq\frac{\sqrt{ 2
\pi} }{ k_n^{1/4} }
\]
and
$\mathbf{1}_{\mathcal{B}_{n}} =
\mathbf{1}_{\mathcal{B}_{n}}(\mathcal{S}_n)$ be its indicator function.
Then
\[
\mathbb{E} \Bigl( \sup_{t > 0} \bigl\llvert
\widetilde{L}_b(t) - J_b(t)\bigr\rrvert
\mathbf{1}_{\mathcal{B}_n} \Bigr) \leq\frac{\sqrt{ 2
\pi} }{ k_n^{1/4} },
\]
and, using (\ref{eqromano}) and the fact that $\sup_t \llvert  \widetilde
{L}_b(t) - J_b(t) \rrvert  \leq1$ almost everywhere,
\[
\mathbb{E} \Bigl( \sup_{t > 0} \bigl\llvert
\widetilde{L}_b(t) - J_b(t)\bigr\rrvert
\mathbf{1}_{\mathcal{B}^c_n} \Bigr) \leq\mathbb{P}\bigl( \mathcal
{B}_n^c \bigr) \leq\frac{1}{k_n^{1/4}}.
\]
Thus
%
\begin{equation}
\label{eqUstat} \mathbb{E} \Bigl( \sup_{t > 0} \bigl\llvert
\widetilde{L}_b(t) - J_b(t) \bigr\rrvert \Bigr) = O
\biggl( \frac{1}{k_n^{1/4}} \biggr) = O \biggl( \frac{b}{n}
\biggr)^{1/4}.
\end{equation}
%

We can now prove the claim of the theorem. First notice that the following
bounds hold:
\begin{eqnarray*}
\mathbb{P}\bigl(\mm{\mathbf{H}}({{\mathcal S}_n},{{\mathbb M}}) >
c_b\bigr) &\leq& \mathbb{P}\bigl(\mm{\mathbf{H}}({{\mathcal
S}_b},{{\mathbb M}}) > c_b\bigr)
\\
&=& J_b(c_b)
\\
&\leq&\widetilde{L}_b(c_b) + \sup_{t > 0}
\bigl\llvert \widetilde{L}_b(t) - J_b(t)\bigr\rrvert.
\end{eqnarray*}
Now take expectations on both sides, and use (\ref{eqLbcb}) and
(\ref{eqUstat}) to obtain that
\[
\mathbb{P}\bigl(\mm{\mathbf{H}}({{\mathcal S}_n},{{\mathbb M}}) >
c_b\bigr) \leq O \biggl( \frac{b}{n} \biggr)^{1/4} +
O \biggl( \frac{1}{n \log n} \biggr) = O \biggl( \frac{b}{n}
\biggr)^{1/4},
\]
as claimed.
\end{pf*}

\begin{pf*}{Proof of Lemma \protect\ref{lemmaconcentration}}
Let $\mathcal{C}= \{c_1, \ldots, c_N\}$ be the set of centers of
Euclidean balls
$\{B_1, \ldots, B_N\}$, forming a minimal $t/2$-covering set for
${\mathbb M}$, for $t/2
< \operatorname{diam}({\mathbb M})$.
Then $\mm{\mathbf{H}}({\mathcal{C}},{{\mathbb M}}) \leq t/2$ and
\begin{eqnarray*}
\mathbb{P} \bigl(\mm{\mathbf{H}}({{\mathcal S}_n},{{\mathbb M}}) > t
\bigr) &\le& \mathbb{P} \bigl(\mm{\mathbf{H}}({{\mathcal S}_n},{
\mathcal{C}}) + \mm{\mathbf{H}}({\mathcal{C}},{{\mathbb M}}) > t \bigr)
\\
&\leq& \mathbb{P} \bigl(\mm{\mathbf{H}}({{\mathcal S}_n},{C}) > t/2
\bigr)
\\
&=& \mathbb{P}(B_j\cap{\mathcal S}_n = \varnothing\mbox{ for some }j)
\\
&\leq&\sum_j P(B_j\cap{\mathcal
S}_n = \varnothing)
\\
&=& \sum_j \bigl[1-P(B_j)
\bigr]^n
\\
&\leq& N \bigl[ 1 - \rho(t) t^d \bigr]^n
\\
&\leq& N \exp \bigl( - n \rho(t) t^d \bigr),
\end{eqnarray*}
where the second-to-last inequality follows from the fact that $\min_j
P(B_j)\geq\rho(t) t^d $ by definition of $\rho(t)$ and the last
inequality from
the fact that $\rho(t) t^d \leq1$.
Next,
let $\mathcal{D} = \{d_1, \ldots, d_{N'} \}$
be the set of centers of $\{B'_1 \cap{\mathbb M}, \ldots, B'_{N'}
\cap{\mathbb M}\}$,
a maximal $t/4$-packing
set
for ${\mathbb M}$. Then $N \leq N'$ [see, e.g., Lemma 5.2,
\citet{niyogi2008finding}, for a proof of this standard fact], and
by definition,
the balls $\{ B'_j \cap{\mathbb M}, j=1,\ldots,N'\}$ are disjoint. Therefore,
\[
1 = P({\mathbb M}) \geq\sum_{j=1}^{N'} P
\bigl(B'_j \cap{\mathbb M}\bigr) = \sum
_{j=1}^{N'} P\bigl(B'_j
\bigr) \geq N' \rho(t/2) \frac{t^d}{2^d}, %
\]
where we have used again the fact that $\min_j P(B'_j) \geq\rho(t/2)
\frac{t^d}{2^d}$. We conclude that $N \leq N' \leq(2^d)/(\rho(t/2) t^d)$.
Hence
\[
\mathbb{P}\bigl(\mm{\mathbf{H}}({{\mathcal S}_n},{{\mathbb M}}) > t
\bigr) \leq \frac
{2^d}{\rho(t/2) t^d}\exp \bigl( -n \rho(t)t^d \bigr).
\]
Now suppose that $t < \min\{ \rho/(2C_2),t_0\}$
($C_2$ and $t_0$ are defined in Assumption~~\ref{assA2}).
Then, since we assume that $\rho(t)$ is differentiable on $[0,t_0]$ with
derivative bounded in absolute value by $C_2$, there exists a $0\leq
\tilde t
\leq t$ such that
\[
\rho(t) = \rho+ t \rho'(\tilde{t}) \geq\rho- C_2t
\geq\frac
{\rho}{2}. %
\]
Similarly, under the same conditions on $t$, we have that $\rho(t/2)
\geq
\frac{\rho}{2}$.
The result follows.
\end{pf*}

\begin{pf*}{Proof of Theorem \protect\ref{theoremahat}}
Let
\[
\hat\rho(x,t) = \frac{P_n(B(x,t/2))}{t^d}. %
\]
Note that
%
\begin{equation}
\label{eqdefC} \sup_{x \in{\mathbb M}} P\bigl(B(x,r_n/2)\bigr)
\leq C r_n^d
\end{equation}
for some $C>0$, since $\rho(x,t)$ is bounded
by Assumption~\ref{assA2}.
Let
$r_n =\break  ( \frac{\log n}{n}  )^{1/(d+2)}$,
and consider all $n$ large enough so that
%
\begin{equation}
\label{eqn} \frac{\rho}{2} ( \log n )^{d/(d+2)} > 1\quad \mbox{and}
\quad n^{2/(d+2)} - \frac{2}{d+2} \log n >0.
\end{equation}
Let $\mathcal{E}_{1,n}$ be the event that the sample $\mathcal{S}_n$
forms an
$r_n$-cover for ${\mathbb M}$. Then, by Lemma~\ref
{lemmaconcentration} and since
$\mathbb{P}(\mathcal{E}^c_{1,n}) = \mathbb{P}(\mm{\mathbf
{H}}({S_n},{{\mathbb M}}) \leq r_n)$, we
have
\begin{eqnarray*}
\mathbb{P}\bigl(\mathcal{E}^c_{1,n}\bigr) &\leq&
\frac{2^{d+1}}{\rho} \biggl( \frac{n}{\log n} \biggr)^{d/(d+2)} \exp \biggl\{
-\frac{\rho}{2} n \biggl( \frac{\log n} {n} \biggr)^{d/(d+2)} \biggr\}
\\
&\leq&\frac{2^{d+1}}{\rho} n^{d/(d+2)} \exp \biggl\{ - \biggl(
\frac{\rho}{2} ( \log n )^{d/(d+2)} \biggr) n^{2/(d+2)} \biggr\}
\\
&\leq&\frac{2^{d+1}}{\rho} n^{d/(d+2)} \exp \bigl\{ - n^{2/(d+2)} \bigr
\}
\\
&\leq&\frac{2^{d+1}}{\rho} \exp \biggl\{ - n^{2/(d+2)} + \frac{d}{d+2}
\log n \biggr\}
\\
&\leq&\frac{2^{d+1}}{\rho} \frac{1}{n},
\end{eqnarray*}
where the third and last inequalities hold since $n$ is assumed large enough
to satisfy~(\ref{eqn}).

Let $C_3 = \max\{C_1,C_2\}$ where
$C_1$ and $C_2$ are defined in Assumption~\ref{assA2}.
Let
\[
\epsilon_n = \sqrt{\overline{C} \frac{2\log n}{(n-1) r_n^d} }
\]
with $\overline{C} = \max\{ 4C_3, 2 (C+1/3)\}$, for some $C$
satisfying~(\ref{eqdefC}).
Assume further that $n$ is large
enough so that, in addition to (\ref{eqn}), $\epsilon_n <1$. With
this choice
of $\epsilon_n$, define the event
\[
\mathcal{E}_{2,n} = \biggl\{ \max_{i=1,\ldots,n} \biggl
\llvert \frac{P_{i,n-1}(B(X_i,r_n/2))}{r_n^d} - \frac{P(B(X_i,r_n/2))}{r_n^d} \biggr\rrvert \leq c^*
\epsilon_n \biggr\},
\]
where $P_{i,n-1}$ is the empirical measure corresponding to the data points
$\mathcal{S}_n \setminus\{X_i\}$, and $c^*$ is a positive number satisfying
\[
c^* \epsilon_n r_n^d - \frac{1}{n}
\geq\epsilon_n r_n^d
\]
for all $n$ large enough (which exists by our choice of $\epsilon_n$ and
$r_n$).
We will show that $\mathbb{P}(\mathcal{E}_{2,n}) \geq1 - \frac{2}{n}$.
To this end, let $\mathbb{P}^i$ denote the probability induced by
$\mathcal{S}_n\setminus\{X_i\}$ (which, by independence is also the
conditional
probability of $\mathcal{S}_n$ given $X_i$) and $\mathbb{P}_i$ be the marginal
probability induced by $X_i$, for $i=1,\ldots,n$.
Then
%
\begin{eqnarray} \label{equnionbound}
\qquad\mathbb{P}\bigl(\mathcal{E}^c_{2,n}\bigr) &\leq&\sum
_{i=1}^n \mathbb {P} \biggl( \biggl\llvert
\frac{P_{n}(B(X_i,r_n/2))}{r_n^d} - \frac{P(B(X_i,r_n/2))}{r_n^d} \biggr\rrvert > c^* \epsilon_n
\biggr)
\nonumber
\\
&\leq&\sum_{i=1}^n \mathbb{P} \biggl(
\bigl\llvert P_{i,n-1}\bigl(B(X_i,r_n/2)\bigr) -
P\bigl(B(X_i,r_n/2)\bigr) \bigr\rrvert > c^*
\epsilon_n r_n^d - \frac{1}{n} \biggr)
\nonumber\\[-8pt]\\[-8pt]
&\leq&\sum_{i=1}^n \mathbb{P} \bigl(
\bigl\llvert P_{i,n-1}\bigl(B(X_i,r_n/2)\bigr) -
P\bigl(B(X_i,r_n/2)\bigr) \bigr\rrvert >
\epsilon_n r_n^d \bigr)
\nonumber
\\
&=& \sum_{i=1}^n \mathbb{E}_i
\bigl[ \mathbb{P}^i \bigl( \bigl\llvert P_{i,n-1}
\bigl(B(X_i,r_n/2)\bigr) - P\bigl(B(X_i,r_n/2)
\bigr) \bigr\rrvert > \epsilon_n r_n^d \bigr)
\bigr],\nonumber
\end{eqnarray}
where the first inequality follows from the union bound, the second
from the
fact
that $\llvert P_{n}(B(X_i,r_n/2)) - P_{i,n-1}(B(X_i,r_n/2))\rrvert  \leq
\frac
{1}{n}$ for all
$i$ (almost\vspace*{1pt} everywhere with respect to the join distribution of the
sample) and
the third inequality from the definition of $c^*$.

By Bernstein's inequality, for any $i =1,\ldots,n$,
\begin{eqnarray*}
&& \mathbb{P}^i \bigl( \bigl\llvert P_{i,n-1}
\bigl(B(X_i,r_n/2)\bigr) - P\bigl(B(X_i,r_n/2)
\bigr) \bigr\rrvert > \epsilon_n r_n^d \bigr)
\\
&&\qquad \leq 2 \exp \biggl\{-\frac{1}{2} \frac{(n -1) \epsilon_n^2 r_n^{2d}}{ C
r_n^d + 1/3
r_n^d \epsilon_n} \biggr\}
\\
&&\qquad \leq 2 \exp \biggl\{ -\frac{1}{2} \frac{(n-1) \epsilon_n^2
r_n^{d}}{ (C +
1/3)} \biggr\}
\\
&&\qquad \leq \frac{2}{n^2},
\end{eqnarray*}
where in the first inequality we have used the fact that $P(B(x,r_n/2))
 ( 1-
P(B(x,r_n/2))  ) \leq C r_n^d$.
Therefore, from (\ref{equnionbound}),
\[
\mathbb{P}\bigl(\mathcal{E}^c_{2,n}\bigr) \leq
\frac{2}{n}.
\]

Let $j=\arg\min_i \hat\rho(X_i,r_n)$ and $k=\arg\min_i \rho
(X_i,r_n)$.
Suppose $\mathcal{E}_{2,n}$ holds and, arguing by contradiction, that
%
\begin{equation}
\bigl\llvert \hat\rho(X_j,r_n) -
\rho(X_k,r_n) \bigr\rrvert = \bigl\llvert \hat {
\rho}_n - \rho(X_k,r_n) \bigr\rrvert > c^*
\epsilon_n.
\end{equation}
Since $\mathcal{E}_{2,n}$ holds, we have
$\llvert \hat\rho(X_j,r_n) - \rho(X_j,r_n) \rrvert  \leq c^* \epsilon
_n $ and
$\llvert \hat\rho(X_k,r_n) - \rho(X_k,r_n) \rrvert  \leq c^* \epsilon_n$.
This implies that if $\hat\rho(X_j,r_n) < \rho(X_k,r_n)$, then
$\rho(X_j,r_n) <
\rho(X_k,r_n) $, while if $\hat\rho(X_j,r_n) > \rho(X_k,r_n)$, then
$\hat\rho(X_k,r_n) < \hat\rho(X_j,r_n) $, which is a contradiction.

Therefore, with probability
$\mathbb{P}  ( \mathcal{E}_{1,n} \cap\mathcal{E}_{2,n}
) \geq
1 - \frac{(2^{d+1}/\rho)+2}{n}$,
the sample points~$\mathcal{S}_n$ forms an $r_n$-covering of ${\mathbb
M}$ and
\[
\Bigl\llvert \hat{\rho}_n - \min_i
\rho(X_i,r_n) \Bigr\rrvert \leq c^* \epsilon_n.
\]
Since the sample $\mathcal{S}_n$ is a covering of ${\mathbb M}$,
%
\begin{equation}
\label{eqmin-rho} \qquad\Bigl\llvert \min_i \rho(X_i,r_n)
- \inf_{x \in{\mathbb M}} \rho(x,r_n) \Bigr\rrvert \leq\max
_i \sup_{x \in B(X_i,r_n)} \bigl\llvert
\rho(x,r_n) - \rho(X_i,r_n)\bigr\rrvert.
\end{equation}
Because $\rho(x,t)$ has a bounded continuous derivative in $t$
uniformly over
$x$, we have, if $r_n < t_0$,
\[
\sup_{x \in B(X_i,r_n)} \bigl\llvert \rho(x,r_n) -
\rho(X_i,r_n)\bigr\rrvert \leq C_3
r_n,
\]
almost surely. Furthermore, since $\rho(t)$ is right-differentiable at zero,
\[
\bigl\llvert \rho(r_n) - \rho\bigr\rrvert \leq C_3
r_n,
\]
for all $r_n < t_0$.
Combining the last two observations with (\ref{eqmin-rho}) and using the
triangle inequality, we conclude that
\[
\llvert \hat{\rho}_n - \rho\rrvert \leq c^*
\epsilon_n + 2C_3 r_n,
\]
with probability at least
$1 - \frac{(2^{d+1}/\rho)+2}{n}$,
for all $n$ large enough.
Because our choice of $r_n$ satisfies the equation
$\epsilon_n = \Theta (\sqrt{ \frac{\log n}{n r_n^d}}  )$,
the terms on
the right-hand side of the
last display are balanced, so that
\[
\llvert \hat{\rho}_n - \rho\rrvert \leq C_4 \biggl(
\frac{\log n}{n} \biggr)^{1/(d+2)},
\]
for some $C_4 >0$.
\end{pf*}

\begin{pf*}{Proof of Theorem \protect\ref{theoremahat2}}
Let $\mathbb{P}_{1}$
denote the unconditional probability measure induced by the first
random half
$\mathcal{S}_{1,n}$ of the sample, $\mathbb{P}_2$ the conditional probability
measure
induced by the second half of the
sample $\mathcal{S}_{2,n}$ given the outcome of the data splitting
and the
values of the first sample, and $\mathbb{P}_{1,2}$ the probability
measure induced by the whole sample and the random splitting. Then
\[
\mathbb{P}\bigl(\mathbf{H}(\mathcal{S}_{2,n},{\mathbb M}) > \hat{t}_{1,n}\bigr) = \mathbb{P}_{1,2}\bigl(\mathbf{H}(
\mathcal{S}_{2,n},{\mathbb M}) > \hat{t}_{1,n}\bigr) =
\mathbb{E}_{1} \bigl( \mathbb{P}_{2}\bigl(\mathbf{H}(
\mathcal {S}_{2,n},{\mathbb M}) > \hat{t}_{1,n}\bigr)
\bigr),
\]
where $\mathbb{E}_1$ denotes the expectation corresponding to $\mathbb{P}_1$.
By Theorem~\ref{theoremahat}, there exist constants $C$ and $C'$
such that the
event
$\mathcal{A}_n$ that $\llvert \hat{\rho}_{1,n} - \rho\rrvert  \leq
C  ( \frac{\log n}{n} )^{1/(d+2)} $ has
$\mathbb{P}_1$-probability no smaller than $ 1- \frac{C'}{n}$, for
$n$ large
enough.
Then
\begin{eqnarray}\label{eqmastereqtheoremhat2}
&& \mathbb{E}_{1} \bigl( \mathbb{P}_{2}\bigl(
\mathbf{H}(\mathcal {S}_{2,n},{\mathbb M}) > \hat{t}_{1,n}
\bigr) \bigr)
\nonumber\\[-8pt]\\[-8pt]\nonumber
&&\qquad \leq \mathbb{E}_{1} \bigl( \mathbb{P}_{2}
\bigl(\mathbf{H}(\mathcal {S}_{2,n},{\mathbb M}) > \hat{t}_{1,n}
\bigr);\mathcal{A}_n \bigr) + \mathbb{P}_1\bigl(
\mathcal{A}_n^c\bigr),\nonumber
\end{eqnarray}
where, for a random variable $X$ with expectation $\mathbb{E}_X$ and
an event
$\mathcal{E}$
measurable with respect to $X$, we write $\mathbb{E}[X;\mathcal{E}]$
for the
expectation of $X$ restricted to the sample points in $\mathcal{E}$.

Define $F(t, \rho) = \frac{2^{d+1}}{t^d \rho} e^{ -((\rho n)/2) t^d }$.
Then, by Lemma~\ref{lemmaconcentration},
\[
\mathbb{P}_2\bigl(\mathbf{H}({\mathcal S}_{2,n},{\mathbb
M}) > \hat {t}_{1,n}\bigr) \leq F(\hat{t}_{1,n},
\rho).
\]
The\vspace*{2pt} rest of the proof is devoted to showing that, on $\mathcal{A}_n$,
$F(\hat{t}_{1,n}, \rho) \leq\alpha+ O  ( (\log n
/2)^{1/(d+2)}  )$.
To simplify the
derivation, we will write $O(R_n)$ to indicate a term that is in
absolute value of order $O  ( (\log n /n)^{1/(d+2)}  ) $,
where the exact value of the constant may change from line to
line. Accordingly, on the event $\mathcal{A}_n$ and for $n$ large
enough, $\hat{\rho}_{1,n} - \rho= O(R_n)$. Since $\rho>0$ by
assumption, this implies that, on the same event and for $n$ large,
\[
\frac{\hat{\rho}_{1,n}}{\rho} = 1 + O(R_n)\quad\mbox{and}\quad \frac{\rho
}{\hat{\rho}_{1,n}} = 1 +
O(R_n). %
\]

Now, on $\mathcal{A}_n$ and for all $n$ large enough,
%
\begin{eqnarray}\label{eqprova}
F(\hat{t}_{1,n}, \rho) &=& \frac{2^{d+1}}{\hat{t}_{1,n}^d \rho} \exp \biggl( -
\frac{n
\hat{t}_{1,n}^d
\rho}{2} \biggr)
\nonumber
\\
&=& \biggl(\frac{\hat{\rho}_{1,n}}{\rho} \biggr) \frac{2^{d+1}}{\hat{t}_{1,n}^d \hat{\rho}_{1,n}} \exp \biggl( -
\frac{n
\hat{t}_{1,n}^d \hat{\rho}_{1,n}
 (\rho/{\hat{\rho}_{1,n}} )}{2} \biggr)
\nonumber
\\
&=& \bigl(1+ O(R_n)\bigr) F(\hat{t}_{1,n},\hat{ \rho}_{1,n}) \exp \biggl( - \frac{n
\hat{t}_{1,n}^d \hat{\rho}_{1,n} O(R_n)}{2} \biggr)
\\
&=& \alpha\bigl(1+ O(R_n)\bigr) \biggl[\exp \biggl( -
\frac{n \hat{t}_{1,n}^d
\hat{\rho}_{1,n} }{2} \biggr) \biggr]^{O(R_n)}
\nonumber
\\
&=& \alpha\bigl(1+ O(R_n)\bigr) \biggl[ \frac{\alpha}{2}
\hat{t}_{1,n}^d \hat{\rho}_{1,n}
\biggr]^{O(R_n)},\nonumber
\end{eqnarray}
where the last two identities follow from the fact that
%
\begin{equation}
\label{eqF} F(\hat{t}_{1,n},\hat{\rho}_{1,n}) =
\alpha
\end{equation}
for all $n$.

Next, let $t^*_n =  (\frac{2}{\alpha\rho} \frac{\log n}{n}
 )^{1/d}$. We then claim that, for all $n$ large enough
and on the
event $\mathcal{A}_n$, $\hat{t}_{1,n} \leq t^*_n$. In fact, using similar
arguments,
\[
F\bigl(t^*_n,\hat{\rho}_{1,n}\bigr) = F
\bigl(t^*_n,\rho\bigr) \bigl(1+O(R_n)\bigr)\exp \bigl( -n
t^*_n O(R_n) \bigr) \rightarrow0\qquad\mbox{as } n \to
\infty,
\]
since $O(R_n) = o(1)$ and, by Lemma~\ref{lemmaconcentration},
$F(t^*_n,\rho) \to0$ and $n \to\infty$.

By (\ref{eqF}), it then follows that, for all $n$ large enough,
$F(\hat{t}_{1,n},\hat{\rho}_{1,n}) >\break F(t^*_n,\hat{\rho}_{1,n})$. Because
$F(t,\rho)$ is decreasing in $t$ for each $\rho$, the claim is proved.

Thus, substituting $\hat{t}_{1,n}$ with $t^*_n$ in equation (\ref
{eqprova})
yields
\begin{eqnarray*}
F(\hat{t}_{1,n}, \rho) &\leq& F\bigl(t^*_n,\rho\bigr)
\\
&\leq& \alpha\bigl(1+ O(R_n)\bigr) \biggl[ \frac{\alpha}{2}
\bigl(t^*_n\bigr)^d \hat{\rho}_{1,n}
\biggr]^{O(R_n)}
\\
&=& \alpha\bigl(1+ O(R_n)\bigr) \biggl[ \frac{\alpha}{2}
\bigl(t^*_n\bigr)^d \rho\bigl(1 + O(R_n)\bigr)
\biggr]^{O(R_n)}
\\
&=& \alpha\bigl(1 + O(R_n)\bigr) \biggl[
\frac{\log n}{n} + o(1) \biggr]^{O(R_n)}
\\
&=& \alpha\bigl(1+O(R_n)\bigr) \bigl(1 + o(1)\bigr)
\\
&=& \alpha+
O(R_n),
\end{eqnarray*}
as $n \to\infty$, where we have written as $o(1)$
all terms that are lower order than $O(R_n)$. The second-to-last step follows
from the limit
\begin{eqnarray*}
\lim_{n\rightarrow\infty} \biggl( \frac{\log n}{n} \biggr)^{O(R_n)}
&=& \lim_{n \rightarrow\infty} \exp \biggl\{ \log \biggl( \frac{\log
n}{n}
\biggr) C \biggl( \frac{\log n}{n} \biggr)^{1/(d+2)} \biggr\}
\\
&=& \exp \biggl\{ \lim_{n \rightarrow\infty} \log \biggl( \frac{\log n}{n}
\biggr) C \biggl( \frac{\log n}{n} \biggr)^{1/(d+2)} \biggr\} = 1,
\end{eqnarray*}
for some constant $C$.

Therefore, on the event $\mathcal{A}_n$ and for all $n$ large enough,
\[
\mathbb{P}_2\bigl(\mathbf{H}({\mathcal S}_{2,n},{\mathbb
M}) > \hat{t}_{1,n}\bigr) \leq\alpha+ O(R_n).
\]
Since $\mathbb{P}_1(\mathcal{A}_n) = O(R_n)$
[in fact, it is of lower order than $O(R_n)$], the result follows from
(\ref{eqmastereqtheoremhat2}).
\end{pf*}

\begin{pf*}{Proof of Theorem \protect\ref{theoremshells}}
Let $B= \sup_{x\in{\mathbb M}}\rho(x,\downarrow0)$.
Fix some $\delta\in(0, B-\rho)$.
Choose equally spaced values
$\rho\equiv\gamma_1 < \gamma_2 < \cdots< \gamma_m < \gamma_{m+1}
\equiv B$
such that
$\delta\geq\gamma_{j+1}-\gamma_j$.
Let $\Omega_j = \{ x\dvtx  \gamma_j \leq\rho(x,\downarrow0) < \gamma
_{j+1}\}$, and
define $  h({\mathcal S}_n, \Omega_j) = \sup_{y \in
\Omega_j} \min_{x
\in{\mathcal
S}_n} \llVert  x-y \rrVert _2$ for $j=1,\ldots, m$.
Now
$
\mm{\mathbf{H}}({{\mathcal S}_n},{{\mathbb M}}) = h({\mathcal
S}_n,{\mathbb M})
\leq\break \max_j h({\mathcal S}_n, \Omega_j)$,
and so
\[
\mathbb{P}\bigl(\mm{\mathbf{H}}({{\mathcal S}_n},{{\mathbb M}}) >t
\bigr) \leq \mathbb{P}\Bigl(\max_j h({\mathcal
S}_n,\Omega_j) >t\Bigr) \leq \sum
_{j=1}^m \mathbb{P}\bigl(h({\mathcal
S}_n,\Omega_j) >t\bigr). %
\]
Let $C_j=\{c_{j1},\ldots, c_{jN_j}\}$ be the set of centers of $\{
B(c_{j1},t/2),\ldots, B(c_{jN_j},t/2)\}$, a $t/2$-covering set for
$\Omega_j$,
for $j=1,\ldots, m$. Let $B_{jk}=B(c_{jk},t/2)$, for $k = 1,\ldots, N_j$.
Then, for all $j$ and for $t < \min\{\gamma_j/(2C_1),t_0\}$
($C_1$ and $t_0$ are defined in Assumption~\ref{assA2}),
there is $0 \leq\tilde{t} \leq t$ such that
\[
\frac{P(B_{jk})}{t^d} = \rho(c_{jk},t) \geq\rho(c_{jk},
\downarrow 0) + t \rho'\bigl(c_{jk}(\tilde{t})\bigr) \geq
\rho(c_{jk}, \downarrow0) - C_1 t \geq \frac{\gamma_j}{2}
\]
and
\begin{eqnarray*}
\mathbb{P}\bigl(h({\mathcal S}_n,\Omega_j) >t\bigr) &\leq&\mathbb {P}\bigl(h({\mathcal S}_n,C_j) +
h(C_j,\Omega_j) >t\bigr)
\\
&\leq&\mathbb{P}\bigl(h({\mathcal
S}_n,C_j) > t/2\bigr)
\\
&\leq&\mathbb{P}(B_{jk} \cap{\mathcal S}_n = \varnothing
\mbox{ for some } k)
\\
&\leq& \sum_{k=1}^{N_j}
\mathbb{P}(B_{jk} \cap{\mathcal S}_n = \varnothing)
\\
&\leq&\sum_{k=1}^{N_j}
\bigl[1-P(B_{jk})\bigr]^n \leq\sum
_{k=1}^{N_j} \biggl[1- t^d
\frac{\gamma_j}{2} \biggr]^n
\\
&=&N_j \biggl[1- t^d \frac{\gamma_j}{2}
\biggr]^n \leq N_j \exp \biggl(- \frac{n
\gamma_j t^d}{2}
\biggr).
\end{eqnarray*}
Following the strategy in the proof of Lemma~\ref{lemmaconcentration},
\[
P(\Omega_j) \geq\sum_k
P(B_{jk})\geq N_j t^d \frac{\gamma_j}{2^{d+1}},
\]
so that $  N_j \leq2^{d+1}
P(\Omega_j)/(\gamma_j t^d)$.
Therefore, for
$t< \min\{\rho/(2C_1),t_0\}$,
\begin{eqnarray*}
\mathbb{P}\bigl(\mm{\mathbf{H}}({{\mathcal S}_n},{{\mathbb M}}) >t
\bigr) &\leq&\frac
{2^{d+1}}{t^d} \sum_{j=1}^m
\frac{P(\Omega_j)}{\gamma_j} \exp \biggl( - \frac{n \gamma_j
t^d}{2} \biggr)
\\
&=& \frac{2^{d+1}}{t^d} \sum_{j=1}^m
\frac{G(\gamma_j+\delta)-G(\gamma_j)}{\delta} \frac{\delta
}{\gamma_j} \exp \biggl( - \frac{n \gamma_j t^d}{2} \biggr)
\\
&\to&\frac{2^{d+1}}{t^d}\int_{\rho}^{B}
\frac{g(v)}{v} \exp \biggl( - \frac{n v
t^d}{2} \biggr)\,dv
\end{eqnarray*}
as $\delta\to0$.
\end{pf*}

\begin{pf*}{Proof of Theorem \protect\ref{theoremghat}}
Let
%
\begin{equation}
g^*(v) = \frac{1}{n}\sum_{i=1}^n
\frac{1}{b} K \biggl( \frac{v - W_i}{b} \biggr),
\end{equation}
where $W_i = \rho(X_i, \downarrow0)$.
Then
\[
\sup_v \bigl\llvert \hat g(v) - g(v)\bigr\rrvert \leq
\sup_v \bigl\llvert g(v) - g^*(v)\bigr\rrvert + \sup
_v \bigl\llvert g^*(v) - \hat g(v)\bigr\rrvert.
\]
By standard asymptotics for kernel estimators,
\[
\sup_v \bigl\llvert g(v) - g^*(v)\bigr\rrvert \leq
C_1 b^2 + O_P \biggl(\sqrt{
\frac{\log n}{n b}} \biggr). %
\]
Next, note that
%
\begin{equation}
\label{eqfindC} \biggl\llvert K \biggl( \frac{v - W_i}{b} \biggr)- K \biggl(
\frac{v -
V_i}{b} \biggr) \biggr\rrvert \leq \frac{C \llvert W_i-V_i\rrvert }{b}\preceq
\frac{r_n}{b}
\end{equation}
from Theorem~\ref{theoremahat}.
Hence,
\[
\bigl\llvert g^*(v) - \hat g(v)\bigr\rrvert \preceq\frac{r_n}{b^2}.
\]
Statement (1) follows.

To show the second statement, we will use the same notation and a similar
strategy as in the proof of Theorem~\ref{theoremahat2}.
Thus
\begin{eqnarray*}
\mathbb{P}\bigl(\mathbf{H}(\mathcal{S}_{2,n},{\mathbb M}) > \hat{t}_{1,n}\bigr) &=& \mathbb{P}_{1,2}\bigl(\mathbf{H}(
\mathcal{S}_{2,n},{\mathbb M}) > \hat{t}_{1,n}\bigr)
\\
&=&
\mathbb{E}_{1} \bigl( \mathbb{P}_{2}\bigl(\mathbf{H}(
\mathcal {S}_{2,n},{\mathbb M}) > \hat{t}_{1,n}\bigr)
\bigr).
\end{eqnarray*}

Let $\mathcal{A}_n$ be the event defined in the proof of Theorem
\ref{theoremahat2} and $\mathcal{B}_n$ the event that $\sup_v
\llvert \hat{g}(v) -
g(v)\rrvert  \leq r_n$.
Then $\mathbb{P}_1(\mathcal{A}_n \cap\mathcal{B}_n) \geq1 -
O(1/n)$. Now,
\begin{eqnarray*}
&& \mathbb{E}_{1} \bigl( \mathbb{P}_{2}\bigl(\mathbf{H}(
\mathcal {S}_{2,n},{\mathbb M}) > \hat{t}_{1,n}\bigr)
\bigr)
\\
&&\qquad  \leq \mathbb{E}_{1} \bigl( \mathbb{P}_{2}\bigl(
\mathbf{H}(\mathcal {S}_{2,n},{\mathbb M}) > \hat{t}_{1,n}
\bigr);\mathcal{A}_n \cap\mathcal{B}_n \bigr) +
\mathbb{P}_1 \bigl( (\mathcal{A}_n \cap
\mathcal{B}_n)^c \bigr).
\end{eqnarray*}
By Theorem~\ref{theoremshells}, conditionally on $\mathcal
{S}_{1,n}$ and the
randomness of the data splitting,
%
\begin{equation}
\label{eqmaster2} \mathbb{P}_{2}\bigl(\mathbf{H}(\mathcal{S}_{2,n},{
\mathbb M}) > \hat{t}_{1,n}\bigr) \leq \frac{2^{d+1}}{\hat{t}_{1,n}^d}\int
_{\rho}^{\infty} \frac
{g(v)}{v} e^{-n v ({\hat{t}_{1,n}^d}/{2})} \,dv.
\end{equation}
We will show next that, on the event $\mathcal{A}_n \cap
\mathcal{B}_n$, the right-hand side of the previous equation is
bounded by $\alpha+ O(r_n)$ as $n \to\infty$. The second claim of
the theorem will then follow. Throughout the proof, we will assume
that the event $\mathcal{A}_n \cap\mathcal{B}_n$ holds.

Recall that $\hat t_{1,n}$ solves the equation
%
\begin{equation}
\label{eqabcde} \frac{2^{d+1}}{\hat{t}_{1,n}^d}\int_{\rho}^{\infty}
\frac
{\hat g(v)}{v} e^{-n v
(\hat{t}_{1,n}^d/2)} \,dv =\alpha
\end{equation}
(and this solution exists for all large $n$).
By assumption, $g(v)$ is uniformly bounded away from 0.
Hence
$g(v)/\hat g(v) = 1+ O(r_n)$ and so
\begin{eqnarray*}
&& \frac{2^{d+1}}{\hat{t}_{1,n}^d}\int_{\rho}^{\infty}
\frac
{g(v)}{v} e^{-n v
(\hat{t}_{1,n}^d/2)} \,dv \nonumber
\\[-2pt]
&&\qquad = \frac{2^{d+1}}{\hat{t}_{1,n}^d}\int
_{\hat{\rho}_n}^{B} \frac{g(v)}{v} e^{-n v
(\hat{t}_{1,n}^d/2)} \,dv +\nonumber
z_n
\\[-2pt]
&&\qquad = \bigl(1+O(r_n)\bigr)\frac{2^{d+1}}{\hat{t}_{1,n}^d}\int_{\hat{\rho
}_n}^{B}
\frac{\hat g(v)}{v} e^{-n v (\hat{t}_{1,n}^d/2)} \,dv + z_n
\\[-2pt]
&&\qquad = \alpha \bigl(1+O(r_n)\bigr) + O(\hat{\rho}_n -
\rho)+ z_n\nonumber
\\[-2pt]
&&\qquad = \alpha+ O(r_n)+ z_n,\nonumber
\end{eqnarray*}
where
\[
z_n = \frac{2^{d+1}}{\hat{t}_{1,n}^d}\int_{\rho}^{\hat{\rho}_n}
\frac{g(v)}{v} e^{-n v
(\hat{t}_{1,n}^d/2)} \,dv. %
\]
Now $\hat t_{1,n} \geq c_1 (\log n/n)^{1/d}\equiv u_n$
for some $c>0$, for all large $n$ since
otherwise the left-hand side of (\ref{eqabcde}) would go to zero, and
(\ref{eqabcde}) would not be satisfied.
So, for some positive $c_2,c_3$,
\begin{eqnarray*}
z_n &\leq& \frac{2^{d+1}}{u_n^d}\int_{\rho}^{\hat{\rho}_n}
\frac
{g(v)}{\rho} e^{-n \rho (u_n^d/2)} \,dv
\leq \llvert \hat\rho_n - \rho\rrvert \frac{c_2}{n^{c_3}} =O
\bigl(\llvert \hat \rho_n - \rho\rrvert \bigr).
\end{eqnarray*}\upqed
\end{pf*}

\begin{pf*}{Proof of Lemma \protect\ref{lemmaHoeff}}
First we show that
\[
\bigl\llvert p_h(x) - p_h(y)\bigr\rrvert \leq
\frac{L }{h^{D+1}}\llVert x-y\rrVert _2. %
\]
This follows from the definition of $p_h$,
\begin{eqnarray*}
&& \bigl\llvert p_h(x) - p_h(y)\bigr\rrvert
\\
&&\qquad = \biggl
\llvert \frac{1}{h^D} \int_{\mathcal X} K \biggl(
\frac
{\llVert  x-u\rrVert  }{h} \biggr) \,dP(u) - \frac{1}{h^D} \int_{\mathcal X}
K \biggl(\frac{\llVert  y-u\rrVert
}{h} \biggr) \,dP(u)\biggr\rrvert
\\
&&\qquad \leq \frac{1}{h^D} \int_{\mathcal X} \biggl\llvert K \biggl(
\frac
{\llVert  x-u\rrVert  }{h} \biggr) - K \biggl(\frac{\llVert  y-u\rrVert  }{h} \biggr) \biggr\rrvert
\,dP(u)
\\
&&\qquad \leq \frac{1}{h^D} \int_{\mathcal X} L \biggl\llvert
\frac{\llVert  x-u\rrVert  }{h} - \frac{\llVert  y-u\rrVert
}{h} \biggr\rrvert \,dP(u)
\\
&&\qquad \leq\frac{L}{h^D} \frac{\llVert  x-y\rrVert  }{h} = \frac{L }{h^{D+1}} \llVert x-y
\rrVert.
\end{eqnarray*}
By a similar argument,
\[
\bigl\llvert \hat p_h(x) - \hat p_h(y)\bigr
\rrvert \leq\frac{L
}{h^{D+1}}\llVert x-y\rrVert. %
\]

Divide
${\mathcal X}$ into a grid
based on cubes
$A_1,\ldots, A_N$
with length of size $\varepsilon$.
Note that $N\asymp(C/\varepsilon)^D$.
Each cube $A_j$ has diameter
$\sqrt{D} \varepsilon$.
Let $c_j$ be the center of $A_j$.
Now
\begin{eqnarray*}
\llVert \hat p_h- p_h\rrVert _{\infty} &=&
\sup_x \bigl\llvert \hat p_h(x) -
p_h(x)\bigr\rrvert = \max_j \sup
_{x\in A_j} \bigl\llvert \hat p_h(x) -
p_h(x)\bigr\rrvert
\\
&\leq& \Bigl(\max_j \bigl
\llvert \hat p_h(c_j) - p_h(c_j)
\bigr\rrvert \Bigr) + 2v,
\end{eqnarray*}
where
$v = \frac{L\varepsilon\sqrt{D} }{2h^{D+1}}$.
We have that
\begin{eqnarray*}
\mathbb{P} \bigl(\llVert \hat p_h - p_h\rrVert
_{\infty} > \delta \bigr) &\leq& \mathbb{P} \Bigl(\max_j
\bigl\llvert \hat p_h(c_j) - p_h(c_j)
\bigr\rrvert > \delta- 2v \Bigr)
\\
&\leq&\sum_j \mathbb{P} \bigl( \bigl\llvert
\hat p_h(c_j) - p_h(c_j)
\bigr\rrvert > \delta- 2v \bigr).
\end{eqnarray*}
Let
$\varepsilon= \frac{\delta h^{D+1}}{2 L \sqrt{D}}$.
Then $2v = \delta/2$, and so
\[
\mathbb{P} \bigl(\llVert \hat p_h - p_h\rrVert
_{\infty} > \delta \bigr) \leq \sum_j
\mathbb{P} \biggl( \bigl\llvert \hat p_h(c_j) -
p_h(c_j)\bigr\rrvert > \frac
{\delta}{2} \biggr).
\]

Note that
$\hat p_h(x)$ is an average of quantities
bounded between zero and
$K(0)/(nh^D)$.
Hence, by Hoeffding's inequality,
\[
\mathbb{P} \biggl( \bigl\llvert \hat p_h(c_j) -
p_h(c_j)\bigr\rrvert > \frac
{\delta
}{2} \biggr) \leq
2 \exp \biggl( - \frac{ n \delta^2 h^{2D}}{2 K^2(0) } \biggr). %
\]
Therefore, summing over $j$, we conclude
\begin{eqnarray*}
\mathbb{P} \bigl(\llVert \hat p_h - p_h\rrVert
_{\infty} > \delta \bigr) &\leq& 2 N \exp \biggl( - \frac{ n \delta^2 h^{2D}}{2 K^2(0) }
\biggr)
\\
&=& 2 \biggl(\frac{C}{\varepsilon} \biggr)^D \exp \biggl( -
\frac{ n
\delta^2 h^{2D}}{2
K^2(0) } \biggr)
\\
&=& 2 \biggl( \frac{4 C L \sqrt{D}}{\delta h^{D+1}} \biggr)^D \exp \biggl( -
\frac{ n \delta^2 h^{2D}}{2 K^2(0) } \biggr).
\end{eqnarray*}\upqed
\end{pf*}

\section{Conclusion}\label{sectionconclusion}

We have presented several methods for
separating noise from signal in persistent homology.
The first three methods are based
on the distance function to the data.
The last uses density estimation.
There is a useful analogy here:
methods that use the distance function to the data
are like statistical methods that use the empirical distribution function.
Methods that use density estimation use smoothing.
The advantage of the former is that it is more directly connected to
the raw data.
The advantage of the latter is that it is less fragile; that is, it is
more robust to noise and outliers.

We conclude by mentioning some open questions that
we plan to address in future work:
\begin{longlist}[(1)]
\item[(1)] We focus on assessing the uncertainty of persistence diagrams.
Similar ideas can be applied to assess uncertainty in barcode plots.
This requires assessing uncertainty at different scales $\varepsilon$.
This suggests examining the variability of
$\mm{\mathbf{H}}({\widehat S_\varepsilon},{S_\varepsilon})$
at different values of $\varepsilon$.
From
\citet{molchanov1998limit},
we have that
%
\begin{equation}
\sqrt{n\varepsilon^D} \mm{\mathbf{H}}({\widehat S_\varepsilon
},{S_\varepsilon}) \rightsquigarrow \inf_{x\in\partial S_\varepsilon} \biggl
\llvert \frac{\mathbb
{G}(x)}{L(x)}\biggr\rrvert,
\end{equation}
where $\mathbb{G}$ is a Gaussian process,
%
\begin{equation}
L(x) = \frac{d}{dt}\inf \bigl\{ p_\varepsilon(y) -p_\varepsilon(x)\dvtx  \llVert x-y\rrVert \leq t \bigr\} \Big|_{t=0}
\end{equation}
and $p_\epsilon(x)$ is the mean of the kernel estimator using a
spherical kernel
with bandwidth $\varepsilon$.
The limiting distribution
it is not helpful for inference because it would be very difficult to estimate
$L(x)$.
We are investigating practical methods for
constructing confidence intervals on
$\mm{\mathbf{H}}({\widehat S_\varepsilon},{S_\varepsilon})$
and using this to assess uncertainty of the barcodes.

\item[(2)] Confidence intervals provide protection
against type I errors, that is, false detections.
It is also important to investigate the power of the methods
to detect real topological features.
Similarly, we would like to quantify the minimax bounds
for persistent homology.

\item[(3)] In the density estimation method, we use a fixed bandwidth.
Spatially adaptive bandwidths might be useful for more refined inferences.

\item[(4)] It is also of interest
to construct confidence intervals for
other topological parameters such as
\emph{the degree $p$ total persistence} defined by
\[
\theta=2\sum d(x,\operatorname{Diag})^p,
\]
where the sum is over the points in ${\mathcal P}$
whose distance from the diagonal is greater than some threshold,
and $\operatorname{Diag}$ denotes the diagonal.

\item[(5)] Our experiments are meant to be a proof of concept.
Detailed simulation studies are needed to see under which
conditions the various methods work well.

\item[(6)] The subsampling method is very conservative due to the fact that
$b=o(n)$.
Essentially, there is a bias of order
$H({\mathcal S}_b,{\mathbb M})-H({\mathcal S}_n,{\mathbb M})$.
We conjecture that it is possible to adjust for this bias.

\item[(7)] The optimal bandwidth for the density estimation method is
an open question.
The usual theory is based on $L_2$ loss which is not necessarily
appropriate for
topological estimation.
\end{longlist}


\section*{Acknowledgments}
The authors would like to thank the anonymous reviewers for their
suggestions and feedback, as well as Fr\'ed\'eric Chazal and Peter Landweber
for their comments.

\begin{supplement}[id=suppA]
\stitle{Supplement to ``Confidence sets for persistence diagrams''\\}
\slink[doi]{10.1214/14-AOS1252SUPP} 
\sdatatype{.pdf}
\sfilename{aos1252\_supp.pdf}
\sdescription{In the supplementary material we give a brief
introduction to persistence homology and provide additional details about homology,
simplicial complexes and stability of persistence diagrams.}
\end{supplement}

%

\printaddresses
\end{document}